\documentclass{article}

\usepackage{arxiv}

\usepackage[utf8]{inputenc} 
\usepackage[T1]{fontenc}    
\usepackage{mathrsfs,amsmath,amssymb,bm}
\usepackage{lipsum}
\usepackage{amsfonts}
\usepackage{epstopdf}
\usepackage{makecell, rotating}
\usepackage{multirow}
\usepackage{multirow}
\usepackage{subfigure}
\usepackage[ruled,linesnumbered]{algorithm2e}
\usepackage{graphicx}
\usepackage[nottoc]{tocbibind}
\usepackage{array}
\usepackage{listings}
\usepackage{cleveref}
\usepackage{caption}
\usepackage{authblk}

\numberwithin{equation}{section}

\title{Solving time-dependent parametric PDEs by multiclass classification-based reduced order model}

\author[1]{Chen Cui\thanks{\texttt{cuichensx@gmail.com}}}
\author[1]{Kai Jiang\thanks{\texttt{kaijiang@xtu.edu.cn}}}
\author[1]{Shi Shu\thanks{\texttt{shushi@xtu.edu.cn}}}

\affil[1]{School of Mathematics and Computational Science, 
Hunan Key Laboratory for Computation and Simulation in Science and Engineering,
Xiangtan University, Xiangtan, Hunan, 411105, China.}

\renewcommand{\headeright}{}
\renewcommand{\undertitle}{}

\begin{document}
\maketitle

\begin{abstract}
  In this paper, we propose a network model, the multiclass classification-based reduced order model (MC-ROM), for solving time-dependent parametric partial differential equations (PPDEs). This work is inspired by the observation of applying the deep learning-based reduced order model (DL-ROM) [14] to solve diffusion-dominant PPDEs. We find that the DL-ROM has a good approximation for some special model parameters, but it cannot approximate the drastic changes of the solution as time evolves. Based on this fact, we classify the dataset according to the magnitude of the solutions and construct corresponding subnets dependent on different types of data. Then we train a classifier to integrate different subnets together to obtain the MC-ROM. When subsets have the same architecture, we can use transfer learning techniques to accelerate offline training. Numerical experiments show that the MC-ROM improves the generalization ability of the DL-ROM both for diffusion- and convection-dominant problems, and maintains the DL-ROM’s advantage of good approximation ability. We also compare the approximation accuracy and computational eﬀiciency of the proper
  orthogonal decomposition (POD) which is not suitable for convection-dominant problems. For diffusion-dominant problems, the
  MC-ROM has better approximation accuracy than the POD in a small dimensionality reduction space, and its computational performance is more efficient
  than the POD's.  
\end{abstract} 
	

\section{Introduction}\label{sec01}

Partial differential equation (PDE) is a fundamental mathematical model in scientific
and engineering computation. It is urgent to develop a numerical approach for solving
PDEs. The approach requires the following properties: high fidelity, generalization
ability (being available to different PDE, different initial-boundary
condition, model parameters, and so on), computational efficiency (being expected to
achieve the optimal $O(N)$ computational complexity). The approach is in terms of
computational PDE model. 
Parametric PDEs (PPDEs) are one of the most important PDEs. Many scientific and
engineering problems, such as control, optimization, inverse design, uncertainty
quantification, Bayesian inference can be described by PPDEs with computational
domains, initial-boundary conditions, source terms, and physical properties as
parameters.  However, numerically solving PPDEs usually requires expensive
computational costs mainly due to multi-query and real-time computing. Therefore,
designing a computational PDE model that meets the above characteristics for the
PPDEs has important applications. However, it is also a
challenge in scientific and engineering computation.

The projection-based linear reduced order model
(ROM)\,\cite{gallivan2002model,hesthaven2014efficient}
is an effective way to improve the computational efficiency of numerically solving
PPDEs. ROM can be divided into offline and online stages. The offline
stage constructs a low-dimensional subspace to approximate the solution manifold using
obtained high fidelity numerical solutions.
The computational tasks on the offline stage are usually expensive.
The online stage obtains an approximated solution
for a new given model parameter based on the low-dimensional subspace.
The proper orthogonal decomposition (POD) method is a popular
algorithm for constructing linear ROM and is effective for many
questions, such as computational fluid dynamics and structural analysis\,\cite{liang2002proper, berkooz1993proper}.
However, the POD method still has some weaknesses, such as 
(i) it requires to construct a relatively high-dimensional subspace to obtain an acceptable numerical solution; (ii) it needs relatively expensive reduction strategies;
and (iii) it has the intrinsic difficulty to handle physical complexity, etc. 
To overcome these difficulties, a non-intrusive and data-driven nonlinear
reduced-order model based on deep learning (or neural network) has been developed\cite{Kutyniok2019,Hesthaven2018}.
 
Using the neural network as an ansatz to solve PDE can be traced back to the late
1990s\,\cite{lagaris1998artificial}. In recent years, with the evolution of the computational power, the explosive development of deep learning has again
attracted much attention of the community of scientific computing. Due to the great
expressivity of neural
network\,\cite{cybenko1989approximation},
the neural-network PDE solver achieves some breakthroughs in solving a single
PDE, especially high-dimensional
PDE\,\cite{raissi2019physics,weinan2018deep,han2018solving}.
Using neural networks to solve PPDEs 
has also been attracted much attention. The idea is to apply
neural networks to learn the parameter to solution mapping. The main works can be divided into two parts based on the steady-state and time-dependent PPDEs.
Firstly, the works on steady-state PPDEs can be divided into supervised
learning\,\cite{Wang2019,
Phillips2020} and unsupervised learning\,\cite{Chen2020,Zhu2019}. Secondly, since
time-dependent PPDEs also involve time variables, the requirements for its generalization ability are also stronger. Solving time-dependent PPDEs using deep learning can be
mainly divided into two perspectives: continuous and discretization. 
Under continuous perspective, the input
and output spaces are both infinite-dimensional, then the corresponding mapping is an
operator between infinite-dimensional spaces.
Refs.\,\cite{benner2015survey,Lu2019,wang2021learning,Li2020c} give
the corresponding algorithms.
Refs.\,\cite{Bhattacharya2020,o2020derivative} introduce dimension reduction technologies to
transform it into a certain finite-dimensional problem. From a discretization perspective, the authors extract the characteristic parameters of the PDE
model, which can be denoted as $\boldsymbol{\mu}=\left[\mu_{1}, \ldots,
\mu_{n_{\mu}}\right] \in \mathcal{P}$, $\mathcal{P}$ is compact in
$\mathbb{R}^{n_{\mu}}$, then design the corresponding network architecture according to
the characteristics of different problems. Refs.\,\cite{Gonzalez2018a,ahmed2020long, chen2018neural} firstly reduce the
dimension of the solution manifold, and then use LSTM or NeuralODE to
evolve low-dimensional system to obtain the solutions at each time layer,
respectively. Ref.\,\cite{berner2020numerically} develops a network model with
spatiotemporal variables as input parameters, which overcomes the curse of dimension.
Fresca et al.\,\cite{Fresca2020} designed a nonlinear ROM based on convolutional
autoencoder, which is called deep learning-based reduced order model (DL-ROM), and use it to
solve some nonlinear time-dependent PDEs. DL-ROM can predict the numerical solution of
the corresponding PPDEs on spatial discretization points for any given $(\boldsymbol{\mu},
t)$. However, the generalization ability of DL-ROM when dealing with
diffusion-dominant problems becomes weak. Designing an effective network that satisfies the properties of the computational PDE model for solving PPDEs is the concern of this article.

In this work, we firstly repeat numerical experiments of one-dimensional ($1 \mathrm{D}$) viscous
Burgers equation in \cite{Fresca2020} by the DL-ROM. We find that when the order of
magnitude of the solution changes drastically as time evolves, DL-ROM cannot capture
the drastic changes of the solution. Based on this observation, we propose 
a multiclass classification-based reduced order model (MC-ROM)
for solving time-dependent PPDEs. MC-ROM classifies the training data according to the
magnitude of the numerical solution and build different subnets for each type of
data. When subsets have the same architecture, we can use transfer learning
techniques to accelerate offline training. Finally, in order to bring a new parameter
into the corresponding subnet in the testing stage, we assign labels to the training
data and train a classifier based on these data by support vector machine (SVM).
MC-ROM is composed of the classifier and these subnets. We also point out that the
choice of subnet depends on the problem. This paper chooses DL-ROM as the subnet
because of its good approximation ability for convection problems and efficient
online calculation. Numerical experiments show that MC-ROM has sufficient generalization
ability for both diffusion- and convection-dominant problems and has the following advantages: (i) MC-ROM can approximate the equation includes some terms that can lead to rapid variation in the solution with acceptable accuracy, while DL-ROM cannot. (ii) MC-ROM is very efficient compared with POD in online computation. We
analyze the online computational complexity of MC-ROM and POD, and both are $O(N)$
for parameter separable problems. However, the latter is much more expensive for
other types of equations. Our numerical results demonstrate that the MC-ROM is
much faster than the POD when solving a two-dimensional (2D) parabolic equation with discontinuous
coefficients. Moreover, POD is not suitable for convection-dominant problems. (iii)
MC-ROM can be easily extended and improved. Once a new model that approximates the
solution mapping is proposed, we can use it as a subnet to improve the approximation
of MC-ROM.

The rest of this paper is organized as follows. Section \ref{sec02} formulates the
problem and gives the corresponding numerical methods. Section \ref{sec03} introduces
the ROM, including POD and  DL-ROM. Section \ref{sec04} finds that when the equation
is diffusion-dominated, DL-ROM's generalization ability has troubles, and we give an
explanation. Based on this observation, in Section \ref{sec05} we propose the MC-ROM
and analyze its online calculation complexity. In Section \ref{sec06}, we apply the
MC-ROM to solve the $1 \mathrm{D}$ viscous Burgers equation and the 2D parabolic
equation. The approximation accuracy and calculation time of solving these equations
are compared with DL-ROM and with POD. Finally, Section \ref{sec07} summarizes and
discusses.

\section{Problem definition}\label{sec02}

Consider the following time-dependent PPDE
\begin{equation}
     \dot{u}(\mathbf{x}, t, \boldsymbol{\mu})+\mathcal{N}[u(\mathbf{x}, t, \boldsymbol{\mu})]+\mathcal{ L}[u(\mathbf{x}, t, \boldsymbol{\mu}); \boldsymbol{\mu}]=0, \quad(\mathbf{x}, t, \boldsymbol{\mu}) \in \Omega \times \mathcal{T} \times \mathcal{P},
     \label{chouxiang}
\end{equation}
where $\dot{u}$ is the derivative of $u$ with respect to $t$, $\Omega \subset \mathbb{R}^{d}$ is a bounded open set, $ \mathcal{T}=[0, T], $ parameter space $ \mathcal{P} \subset \mathbb{R}^{n_{\mu}}$ is compact, and $\mathcal{N}, \mathcal{L}$ are nonlinear and linear operator respectively. This equation\,\eqref{chouxiang} is well-defined with appropriate initial and boundary conditions.

Use spatial discretization methods such as finite element method (FEM) and finite difference method (FDM) to discretize equations \eqref{chouxiang}, we obtain the following ODE system
\begin{equation}
     \dot{\mathbf{u}}_{h}(t, \boldsymbol{\mu})= \boldsymbol{f}(\mathbf{u}_{h}, t; \boldsymbol{\mu})\quad(t, \boldsymbol{\mu}) \in \mathcal{T} \times \mathcal{P},
     \label{chouxiangODE}
\end{equation}
where $\mathbf{u}_{h}: \mathcal{T} \times \mathcal{P} \rightarrow \mathbb{R}^{N_{h}}$, $N_{h}$ is the spatial degrees of freedom.

Next, for time discretization, $ \mathcal{T}$ is uniformly divided into $N_t$ segments, $\left\{t^{k}\right\}_{k=1}^{N_{t}}$ are time layers to be solved, where $ t^{k }=k \Delta t$, $\Delta t=T / N_{t}$. 
We apply linear multi-step method to solve \eqref{chouxiangODE}, which becomes a fully discretized system
\begin{equation}
    \boldsymbol{r}^{k}\left(\mathbf{u}_{h}^{k}; \boldsymbol{\mu}\right)=\mathbf{0}, \quad k=1, \ldots , N_{t},
    \label{FOM}
\end{equation}
the discrete residual of the $k$th time layer $\boldsymbol{r}^{k}: \mathbb{R}^{N_h} \times \mathcal{P} \rightarrow \mathbb{R}^{N_h} $ is defined as
\begin{equation}
    \boldsymbol{r}^{k}:(\boldsymbol{\xi}; \boldsymbol{\mu}) \mapsto \alpha_{0} \boldsymbol{\xi}-\Delta t \beta_{0} \boldsymbol{ f}\left(\boldsymbol{\xi}, t^{k}; \boldsymbol{\mu}\right)+\sum_{j=1}^{K} \alpha_{j} \mathbf{u}_ {h}^{k-j}-\Delta t \sum_{j=1}^{K} \beta_{j} \boldsymbol{f}\left(\mathbf{u}_{h}^{k-j}, t ^{k-j}; \boldsymbol{\mu}\right).
    \label{lmm}
\end{equation}
A concrete linear multi-step method is determined by the choose of coefficients
$\alpha_{j},\,\,\beta_{j}, j=0, \ldots, K$ with $\sum_{j=0}^{K} \alpha_{j}=0$.
Expressions\,\eqref{chouxiangODE}-\eqref{lmm} are full order model (FOM) of
PDE\,\eqref{chouxiang}.

Given a parameter $\boldsymbol{\mu} \in \mathcal{P}$, its corresponding
PDE\,\eqref{chouxiang} is solved by FOM. When the spatial discrete degree of freedom
$N_{h}$ is large, the obtained algebraic system is large-scale. Solving such a large
system is expensive. Moreover, when $\mu$ changes, the FOM is required to repeat the
above expressions\,\eqref{chouxiangODE}-\eqref{lmm}, which increases the
computational cost again. How to solve these problems is the focal point of this paper.  
\section{Reduced order model}\label{sec03}

From the previous section we saw when $N_{h}$ is large, the equation \eqref{FOM} is large-scale. However, the inherent dimension of solution manifold $\mathcal{S}_{h}=\left\{ \mathbf{u}_{h}(t; \boldsymbol{\mu}) \mid t \in[0, T]\right.$ and $\left.\boldsymbol{\mu} \in \mathcal{P } \subset \mathbb{R}^{n_{\mu}}\right\}$ usually not exceed the dimension of parameter space\,\cite{Quarteroni2015,hesthaven2016certified}. Hence we can use the solution on the low-dimensional manifold to accurately represent the FOM solution $\mathbf{u}_{h} (t; \boldsymbol{\mu})$. Based on this hypothesis, ROM was proposed. In this section we introduce a commonly used linear ROM and a nonlinear ROM.

\subsection{Projection based ROM}\label{POD}

Let $n\ll N_h$, $\mathbf{u}_{n}(t; \boldsymbol{\mu}) \in \mathbb{R}^n$ be the dimension and numerical solution of low-dimensional space, respectively. $V \in \mathbb{R}^{N_{h} \times n}$ is a linear lifting operator, $\operatorname{Col}(V)$ represents the linear space spanned by the column vector of $V$. Linear ROM pursuits an approximation on a $n$-dimensional trial manifold $\mathcal{S}_n=\operatorname{Col}(V)$, the approximation of FOM solution $\mathbf{u}_{h}$ can be obtained by a decoder $V$, denoted as $\boldsymbol{h}_{\mathrm{dec}}$
\begin{equation}
    \mathbf{u}_{h}(t; \boldsymbol{\mu}) \approx \tilde{\mathbf{u}}_{h}(t; \boldsymbol{\mu}) =V \mathbf{u }_{n}(t; \boldsymbol{\mu}).
\end{equation}

POD is one of the most widely used method to generate such a linear trial manifold. It corresponds to the principal component analysis (PCA) in data science. The difference between POD and PCA is that the former retains the reduced system and the latter does not. In the following, we will briefly introduce POD. More details can refer to \cite{liang2002proper,Quarteroni2015,hesthaven2016certified}.

\textbf{Step 1: Data collection}. Sample $N_{train}$ parameters in the parameter space in a certain distribution, numerically solving the corresponding PDE\,\eqref{FOM}, and assemble them into a snapshot matrix $S \in \mathbb{R}^{N_{h} \times N_{s}}$
\begin{equation}
     S=\left[\mathbf{u}_{h}\left(t^{1}; \boldsymbol{\mu}_{i}\right) \ldots \mathbf{u}_{h}\left( t^{N_{t}}; \boldsymbol{\mu}_{i}\right), i=1 \ldots N_{\text {train }}\right],
\end{equation}
here $N_{s}=N_{\text {train}} N_{t}$.

\textbf{Step 2: Computing singular value decomposition (SVD) of $S$}
\begin{equation}
     S=U \Sigma Z^{T},
\end{equation}
where $U=\left[\zeta_{1},\ldots, \zeta_{N_{h}}\right] \in \mathbb{R}^{N_{h} \times
N_{h}}, \, Z= \left[\psi_{1},\ldots, \psi_{N_{s}}\right] \in \mathbb{R}^{N_{s} \times
N_{s}}$ are both orthogonal, $ \Sigma=\operatorname{diag}\left(\sigma_{1}, \ldots,
\sigma_{r}\right) \in \mathbb{R}^{N_{h} \times N_{s}}$, singular values $\sigma_{1}
\geq \sigma_{2} \geq \ldots \geq \sigma_{r}>0$, $r \leq \min \left(N_{h}, N_{s}\right
)$. Take the first $n$ left singular vectors of $U$ as $V=\left[
\zeta_{1},\ldots, \zeta_{n}\right]$.

\textbf{Step 3: Generate reduced system.} Substituting $V \mathbf{u}_{n}(t; \boldsymbol{\mu})$ into the ODE system \eqref{chouxiangODE} we get
\begin{equation}
    V \dot{\mathbf{u}}_{n}(t; \boldsymbol{\mu})=\mathbf{f}\left(V \mathbf{u}_{n}(t; \boldsymbol{\mu}), t; \boldsymbol{\mu}\right)\quad(t, \boldsymbol{\mu}) \in \mathcal{T} \times \mathcal{P}.
\end{equation}
Let the residuals be orthogonal to $\operatorname{Col}(V)$, we obtain the reduced system
\begin{equation}
    V^{T} V \dot{\mathbf{u}}_{n}(t; \boldsymbol{\mu})=V^{T} \mathbf{f}\left(V \mathbf{u}_ {n}(t; \boldsymbol{\mu}), t; \boldsymbol{\mu}\right)\quad(t, \boldsymbol{\mu}) \in \mathcal{T} \times \mathcal{P }.
    \label{ROM}
\end{equation}

The POD method introduced above is a Galerkin method. More generally, let $W \in \mathbb{R}^{N_{h} \times n}\,, \,W\neq V$ and
\begin{equation}
    W^{T} V \dot{\mathbf{u}}_{n}(t; \boldsymbol{\mu})=W^{T} \mathbf{f}\left(V \mathbf{u}_ {n}(t; \boldsymbol{\mu}), t; \boldsymbol{\mu}\right)\quad(t, \boldsymbol{\mu}) \in \mathcal{T} \times \mathcal{P },
    \label{ROMPG}
\end{equation}
we obtain another reduced system, which is a Petrov-Galerkin method.

Suppose $U_{h}$ is a Hilbert space spanned by all $N_h$-dimensional vectors with respect to the appropriate inner product. $\mathcal{S}_{h}\subset U_{h}$ is low-dimensional solution manifold. We try to find a $n$-dimensional subspace $\mathcal{S}_{n}$ to approximate $\mathcal{S}_{h}$, the degree of approximation is characterized by the Kolmogorov $n$-width
\begin{equation}\nonumber
    d_{n}\left(\mathcal{S}_{h}; U_{h}\right):=\inf _{\mathcal{S}_{n} \subset U_{h} \atop \operatorname{ dim}\left(\mathcal{S}_{n} \right)=n} d\left(\mathcal{S}_{h}; \mathcal{S}_{n}\right)=\inf _ {\mathcal{S}_{n} \subset U_{h} \atop \operatorname{dim}\left(\mathcal{S}_{n} \right)=n} \sup _{\boldsymbol{\mu }\in\mathcal{P}} \inf _{\mathbf{u}_{n} \in \mathcal{S}_{n}}\left\|\mathbf{u}_{h}(\boldsymbol {\mu})-\mathbf{u}_{n}\right\|_{U_h},
\end{equation}
where $d\left(\mathcal{S}_{h}; \mathcal{S}_{n}\right)$ is the maximum distance between $\mathcal{S}_{n}$ and $\mathcal{S}_{h}$. Theories show that Kolmogorov $n$-width decreases with the increase of $n$, but the decay rate depends on the problem \cite{Ohlberger2015}. For the linear trial subspace $\mathcal{S}_{n}$ generated by POD, the Kolmogorov $n$-width defined above decays quickly for diffusion-dominant problems. However, Kolmogorov $n$-width decays slowly for convection-dominant problems. We must increase the order of ROM largely to obtain the required accuracy\,\cite{Ohlberger2015}. Next, we introduce a convolutional autoencoder based nonlinear ROM.

\subsection{Convolutional autoencoder based nonlinear ROM}
Autoencoder is a kind of feedforward neural network that tries to learn the identical
mapping, which consists of two parts: encoder $\boldsymbol{h}_{\mathrm{enc}}$ and
decoder $\boldsymbol{h}_{\mathrm{dec}}$. 
Encoder maps $\mathbf{u}_{h}$ to a low-dimensional vector $\mathbf{u}_{n}$
\begin{equation}
    \begin{aligned}
        \boldsymbol{h}_{\mathrm{enc}}: \mathbb{R}^{N_{h}} & \rightarrow \mathbb{R}^{n} \\
        \mathbf{u}_{h} & \mapsto \mathbf{u}_{n},
    \end{aligned}
\end{equation}
where $n \ll N_h$.
Decoder maps $\mathbf{u}_n$ to an approximation of FOM solution $\mathbf{u}_{h}$
\begin{equation}
    \begin{aligned}
        \boldsymbol{h}_{\mathrm{dec}}: \mathbb{R}^{n} & \rightarrow \mathbb{R}^{N_{h}} \\
        \mathbf{u}_{n} & \mapsto \widetilde{\mathbf{u}}_{h}.
    \end{aligned}
\end{equation}
Therefore, autoencoder can be written in the following form
\begin{equation}
\boldsymbol{h}: \mathbf{u}_{h} \mapsto \boldsymbol{h}_{\mathrm{dec}} \circ \boldsymbol{h}_{\mathrm{enc}}(\mathbf{u }_{h}).
\end{equation}
The encoder contains $L$ layers network
\begin{equation}
     \left\{\begin{array}{l}
     \boldsymbol{\mathbf{u}}_{n}^{(0)}=\mathbf{u}_{h} \in \mathbb{R}^{N_h} \\
     \boldsymbol{\mathbf{u}}_{n}^{(l)}=\varphi\left(\boldsymbol{h}^l_{\mathrm{enc}}(\boldsymbol{\mathbf{u}}_ {n}^{(l-1)})\right), \quad l=1, \cdots, L-1, \\
     \boldsymbol{\mathbf{u}}_{n}=\boldsymbol{\mathbf{u}_{n}}^{(L)}=\boldsymbol{h}^{L}_{\mathrm{enc} }(\boldsymbol{\mathbf{u}}_{n}^{(L-1)}) \in \mathbb{R}^{n},
     \end{array}\right. 
     \label{enc}
\end{equation}
where $\boldsymbol{\mathbf{u}}_{n}^{(l)}$ represents the $l$-layer output. $\boldsymbol{h}^{l}_{\mathrm{enc}}$ represents the $l$-layer network operator, such as convolution or fully connection. $\varphi $ is a nonlinear activation function.
The decoder contains $L$ layers network
\begin{equation}
     \left\{\begin{array}{l}
     \tilde{\mathbf{u}}_{h}^{(0)}=\boldsymbol{\mathbf{u}_{n}} \in \mathbb{R}^{n} \\
     \tilde{\mathbf{u}}_{h}^{(l)}=\varphi\left(\boldsymbol{h}^l_{\mathrm{dec}}(\tilde{\mathbf{u}}_ {h}^{(l-1)})\right), \quad l=1, \cdots, L-1, \\
     \tilde{\mathbf{u}}_{h}=\tilde{\mathbf{u}}_{h}^{(L)}=\boldsymbol{h}^{L}_{\mathrm{dec} }(\tilde{\mathbf{u}}_{h}^{(L-1)}) \in \mathbb{R}^{N_h}.
     \end{array}\right. 
     \label{dec}
\end{equation}

Different architecture in the autoencoder leads to different method. Fully connected
layers have good expressivity. Refs. \cite{wang2016auto,hartman2017deep,otto2019linearly} choose all layers as fully
connected layers then get multilayer perceptron autoencoder (MLPAE), which is
effective for small-scale systems. However, for large-scale systems, the required
training parameters in MLPAE will increase dramatically, which results in a curse of
dimension. For example, when $N_h$ is about $10^{6}$, even if only one fully
connected layer is used to reduce the high fidelity vector to $100$-dimension, the
network parameters will exceed $10^{8}$. Therefore, it requires to reduce parameters
in MLPAE. An alternative is the convolutional autoencoder (CAE). 
Convolutional networks not only have the characteristics, such as
the connectivity, the translation invariance and so
on \cite{lecun2015deep,goodfellow2016deep}, but also has well approximation
properties\cite{bao2014approximation,zhou2020universality,petersen2020equivalence,he2021approximation}.
Refs. \cite{Zhu2018,Zhu2019} use a pure convolutional network to construct the
ROM which brings the benefit of reducing the data required for network
training. There are also some works to use both convolutional layers and
fully connected layers to construct a ROM\,\cite{Gonzalez2018a, Lee2020,
Fresca2020} which have fewer parameters than MLPAE and are easier to train
than pure convolutional network. In the further work, we will incorporate
these ROM into our MC-ROM models according to considered problems.

\begin{figure}[htb]
    \centering
    \includegraphics[width=1.0\textwidth]{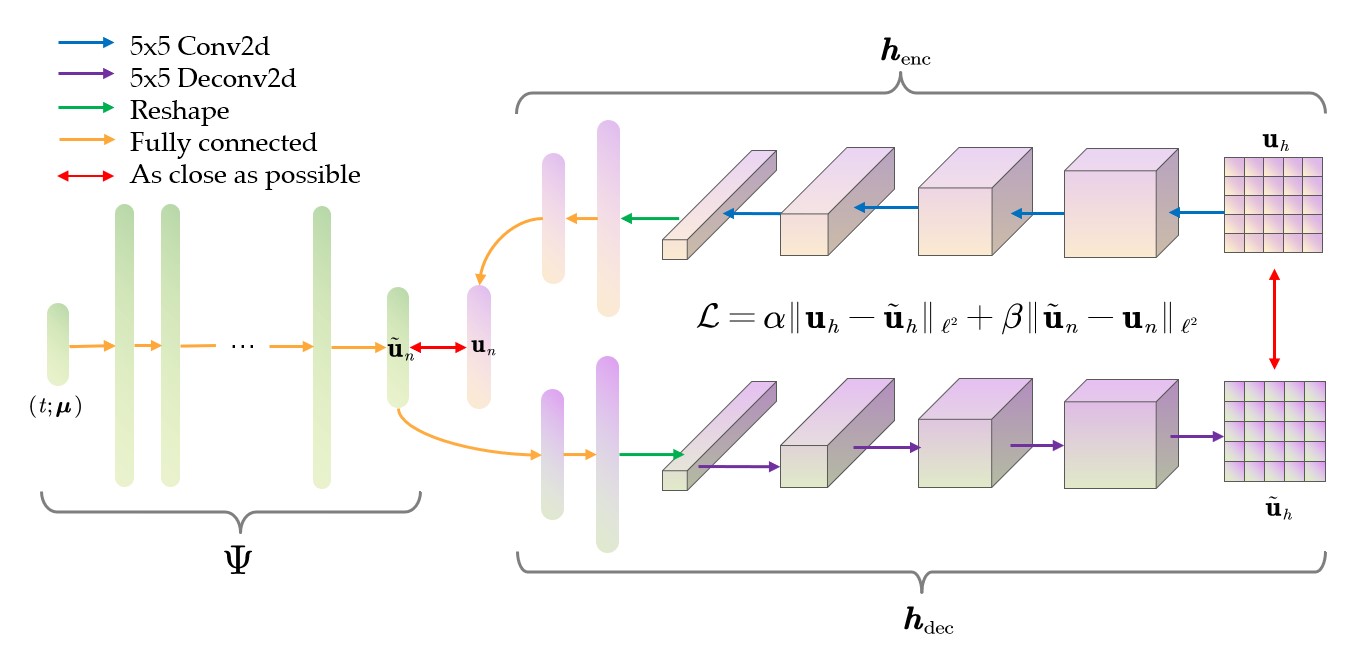}
    \caption{\textbf{DL-ROM}: The network architecture proposed in \cite{Fresca2020}.  \textit{Upper right:} $\boldsymbol{h}_{\mathrm{enc}}$ converts original data $\mathbf{u}_{h}$ into a latent representation $\mathbf{u}_{n}$;  \textit{Lower right:} $\boldsymbol{h}_{\mathrm{dec}}$ recovers data from low-dimensional solution $\tilde{\mathbf{u}}_{n}$; \textit{Left:}$\Psi $ is to fit the low-dimensional solution map.}
    \label{fig:DL_ROM}
\end{figure}

Figure \ref{fig:DL_ROM} shows a ROM based on CAE, called DL-ROM. $\boldsymbol{h}_{\mathrm{enc}}$ represents the encoder. It is a $6$-layers network, where the first four layers are convolutional layers, and the last two layers are fully connected layers. $\boldsymbol{\theta}_E$ represents the network parameters of $\boldsymbol{h}_{\mathrm{enc}}$, then formula \eqref{enc} can be written as
\begin{equation}
    \mathbf{u}_{n}=\boldsymbol{h}_{\operatorname{enc}}\left(\mathbf{u}_{h};\boldsymbol{\theta}_E\right).
    \label{encoder}
\end{equation}$\Psi $ is the low-dimensional solution map
\begin{equation}
    \tilde{\mathbf{u}}_{n}=\Psi(t, \boldsymbol{\mu};\boldsymbol{\theta}_F),
    \label{psi}
\end{equation}
where $\boldsymbol{\theta}_F$ represents the network parameters of $\Psi$. We expect its output $\tilde{\mathbf{u}}_{n}$ be as close to the encoder output $\mathbf{u}_{n }$ as possible. DL-ROM uses $\Phi:=\boldsymbol{h}_{\mathrm{dec}}\circ \Psi $ to approximate the solution map
\begin{equation}
    \tilde{\mathbf{u}}_{h}=\boldsymbol{h}_{\operatorname{dec}}\left(\tilde{\mathbf{u}}_{n};\boldsymbol{\theta} _D\right),
    \label{decoder}
\end{equation}
where $\boldsymbol{\theta}_D$ represents the network parameters of $\boldsymbol{h}_{\mathrm{dec}}$. It is also a $6$-layers network, where the first two layers are fully connected, the last four layers are deconvolutional layers.

The training of DL-ROM is carried out on its three subnets $\boldsymbol{h}_{\mathrm{enc}},$$\boldsymbol{h}_{\mathrm{dec}},$ $\Psi$ simultaneously, hence the loss function includes both the autoencoder error and the low-dimensional fitting error
\begin{equation}
    \min _{\boldsymbol{\theta}} \mathcal{J}(\boldsymbol{\theta})=\min _{\boldsymbol{\theta}} \frac{1}{N_{s}} \sum_{ i=1}^{N_{\text {train }}} \sum_{k=1}^{N_{t}} \mathcal{L}\left(t^{k}, \boldsymbol{\mu}_ {i}; \boldsymbol{\theta}\right),
    \label{error_eq}
\end{equation}where
\begin{equation}
    \mathcal{L}\left(t^{k}, \boldsymbol{\mu}_{i}; \boldsymbol{\theta}\right)=\alpha \left\|\mathbf{u}_{h} \left(t^{k}; \boldsymbol{\mu}_{i}\right)-\tilde{\mathbf{u}}_{h}\left(t^{k}; \boldsymbol{\mu }_{i}, \boldsymbol{\theta}_{F}, \boldsymbol{\theta}_{D}\right)\right\|_{\ell^2}^2+\beta \left\|\tilde {\mathbf{u}}_{n}\left(t^{k}; \boldsymbol{\mu}_{i}, \boldsymbol{\theta}_{F}\right)-\mathbf{u} _{n}\left(t^{k}; \boldsymbol{\mu}_{i}, \boldsymbol{\theta}_{E}\right)\right\|_{\ell^2}^2.
    \label{eq:316}
\end{equation}
$\boldsymbol{\theta}=\left(\boldsymbol{\theta}_{E}, \boldsymbol{\theta}_{D}, \boldsymbol{\theta}_{F}\right)$, $ \alpha,\,\beta \in(0,1)$ are hyperparameters. Algorithm \ref{algorithm:DL_ROM_offline} and Algorithm \ref{algorithm:DL_ROM_online} summarize the training and testing process of DL-ROM, respectively.

\begin{algorithm}
    \caption{DL-ROM offline traning}\label{algorithm:DL_ROM_offline}
    \KwData{Parameter matrix $P \in \mathbb{R}^{\left(n_{\mu}+1\right) \times N_{s}},$ snapshot matrix $S \in \mathbb{R}^{N_{h} \times N_{s}}$}
    \KwIn{The number of training epochs $N_{\text {epochs}},$ batch size $N_{b}$, learning rate $\eta$, learning rate adjustment strategy, early stopping patiences, hyperparameter $\alpha, \beta $} 
    \KwOut{Trained model $\boldsymbol{h}_{\mathrm{enc}},\boldsymbol{h}_{\mathrm{dec}},\Psi$} 
    Data preprocess on $S$ and scale it to $[0,1]$

    Randomly shuffle $P, S$, and split them into training set $\left(P^{\text {train }}, S^{\text {train }}\right)$ and validation set $\left(P^{\text {val}}, S^{\text {val}}\right)$ according to the ratio of 8:2

     Divide data into mini-batches $\left(P^{\text {batch }}, S^{\text {batch }}\right)$ and get data loader

    Define model $\boldsymbol{h}_{\mathrm{enc}},\boldsymbol{h}_{\mathrm{dec}},\Psi$ and initialize the model parameters through the Kaiming
    uniform initialization\cite{he2015delving}

    \While{$n_{\text {epochs}} \leq N_{\text {epochs}}$}{
        \For{$\left(P^{\text {batch }}, S^{\text {batch }}\right)$ in train-loader}{
            Encode $S^{\text {batch }}$ according to equation \eqref{encoder} to get the corresponding low-dimensional representation $S^{\text {batch }}_n$

            Compute equation \eqref{psi} for all $t, \boldsymbol{\mu}$ in $P^{\text {batch}}$ to get an approximation $\widetilde{S^{\text {batch }}_n}$ to $S^{\text {batch }}_n$

            Decode $\widetilde{S^{\text {batch }}_n}$ according to equation \eqref{decoder} to get $\widetilde{S^{\text {batch }}}$

            Computing loss function \eqref{error_eq}

            Backpropagate

            Apply Adam algorithm to update model parameters
        }
        \For{$\left(P^{\text {batch }}, S^{\text {batch }}\right)$ in val-loader}{
            Do 7-10 similar to what was done on training set
            
            \If{early-stopping}{
                \Return{trained model $\boldsymbol{h}_{\mathrm{enc}},\boldsymbol{h}_{\mathrm{dec}},\Psi$}
            } 
        }
        Adjust the learning rate

        $n_{\text {epochs}}=n_{\text {epochs}}+1$.
    }
    \Return{trained model $\boldsymbol{h}_{\mathrm{enc}},\boldsymbol{h}_{\mathrm{dec}},\Psi$}
\end{algorithm}

\begin{algorithm}
    \caption{DL-ROM online testing}\label{algorithm:DL_ROM_online}
    \KwData{Testing parameter matrix $P \in \mathbb{R}^{\left(n_{\mu}+1\right) \times N_{s}}$}
    \KwIn{Trained model $\Psi, \boldsymbol{h}_{\mathrm{dec}}$ by Algorithm \ref{algorithm:DL_ROM_offline}}
    \KwOut{Predicted solution matrix $\widetilde{S}$ corresponding to the parameter matrix $P$}
    Compute equation \eqref{psi} for all $t, \boldsymbol{\mu}$ in $P$ to get $\widetilde{S_n}$

    Decode $\widetilde{S_n}$ according to equation \eqref{decoder} to get $\widetilde{S}$.
\end{algorithm}

\section{Observation}\label{sec04}

In this section, we firstly repeat numerical experiments of the $1 \mathrm{D}$
viscous Burgers equation in \cite{Fresca2020}. Then we change the parameter range and
use DL-ROM to solve the Burgers equation again. Numerical results show that DL-ROM
cannot capture the drastic change of solution when the equation becomes
diffusion-dominant. We give an explanation for this phenomenon.

Consider the following $1 \mathrm{D}$ viscous Burgers equation
\begin{equation}
    \left\{\begin{array}{ll}
    \frac{\partial u}{\partial t}+u \frac{\partial u}{\partial x}-\frac{1}{\mu} \frac{\partial^{2} u}{\partial x^{2}}=0, & (x, t) \in(0, L) \times(0, T) \\
    u(0, t)=0, & t \in(0, T) \\
    u(L, t)=0, & t \in(0, T) \\
    u(x, 0)=u_{0}(x), & x \in(0, L)
    \end{array}\right. 
    \label{burgers}
\end{equation}
where initial value
\begin{equation}
    u_{0}(x)=\frac{x}{1+\sqrt{1 / A_{0}} \exp \left(\mu x^{2} / 4\right)},
\end{equation}
and $A_{0}=\exp (\mu / 8),\, L=1$, $T=2$, $\mu \in \mathcal{P} \subset \mathbb{R}^{n_{\mu }}$ is the single parameter ($n_{\mu}=1$). The numerical methods used here are linear FEM with $N_{h}=256$ grid points, and backward Euler with $N_{t}=100$ time layers.

As Ref \cite{Fresca2020} does, the parameter range is $\mathcal{P}=[100,1000]$, and the corresponding viscosity coefficient $1/\mu$ belongs to $ [10^{-3},10^{-2}]$. This means that equation \eqref{burgers} is a convection-dominated equation. We randomly sample $N_{\text {train}}=20$ parameters in the parameter range according to the uniform distribution. The midpoint of every two training parameters is used as test parameter, i.e., $N_{\text {test}}=19$. After sampling, we solve the corresponding equation \eqref{burgers} by the above numerical methods, and assemble the parameters and corresponding high-fidelity solutions into a parameter matrix $P \in \mathbb{R}^{\left(n_{\mu}+1\right) \times N_{s}},$ and a snapshot matrix $S \in \mathbb{R}^{N_{h} \times N_{s}}$.

We use \lstinline{Pytorch}\cite{paszke2019pytorch} as a implementation platform for training the network and predicting results. The concrete network architecture is as follows. $\Psi$ contains $10$ hidden layers, and each layer contains $50$ neurons. The nonlinear activation function used here is ELU \cite{clevert2015fast}. Table \ref{tab:Network_structure_burgers} shows the parameters of $\boldsymbol{h}_{\mathrm{enc}},\,\boldsymbol{h}_{\mathrm{dec}}$. 
\begin{table}[htb]
    \centering
    \caption{Network architecture used in $1 \mathrm{D}$ viscous Burgers equation experiment. $N_{b}$ represents batch size. Notations in Conv2d and ConvTranspose2d are: i: in channels; o: out channels; k: kernel size; s: stride; p: padding.}
    \resizebox{\textwidth}{!}{
    \begin{tabular}{|llr|ll|}
        \hline
        \multicolumn{2}{|c}{$\boldsymbol{h}_{\mathrm{enc}}$ (input shape:($N_b$, 256))} &        & \multicolumn{2}{c|}{$\boldsymbol{h}_{\mathrm{dec}}$ (input shape:($N_b$, $n$))} \\
        \hline
        Layer type & Output shape &        & Layer type & Output shape \\
        \hline
        Reshape & ($N_b$,1,16,16) &        & Fully connected & ($N_b$,256) \\
        Conv2d(i=1,o=8,k=5,s=1,p=2) & ($N_b$,8,16,16) &        & Fully connected & ($N_b$,256) \\
        Conv2d(i=8,o=16,k=5,s=2,p=2) & ($N_b$,16,8,8) &        & Reshape & ($N_b$,64,2,2) \\
        Conv2d(i=16,o=32,k=5,s=2,p=2) & ($N_b$,32,4,4) &        & ConvTranspose2d(i=64,o=64,k=5,s=3,p=2) & ($N_b$,64,4,4) \\
        Conv2d(i=32,o=64,k=5,s=2,p=2) & ($N_b$,64,2,2) &        & ConvTranspose2d(i=64,o=32,k=5,s=3,p=1) & ($N_b$,32,12,12) \\
        Reshape & ($N_b$,256) &        & ConvTranspose2d(i=32,o=16,k=5,s=1,p=1) & ($N_b$,16,14,14) \\
        Fully connected & ($N_b$,256) &        & ConvTranspose2d(i=16,o=1,k=5,s=1,p=1) & ($N_b$,1,16,16) \\
        Fully connected & ($N_b$,$n$)  &        & Reshape & ($N_b$,256) \\
        \hline
    \end{tabular}%
    }
    \label{tab:Network_structure_burgers}%
\end{table}%

We use Adam algorithm \cite{kingma2014adam} to train the network, the initial
learning rate is set to $\eta = 10^{-4}$, the batch size is $N_{b}=20$. We set early
stopping patiences to $500$. In other words, we stop training if error on the
validation set does not decrease for $500$ consecutive epochs. All these
settings are the same as DL-ROM in \cite{Fresca2020}. Under these settings, we apply Algorithm \ref{algorithm:DL_ROM_offline} for training. In the case of a fixed time-parameter instance, define its relative error as
\begin{equation}
    \mathcal{E}_{single}= \frac{\left\|\mathbf{u}_{h}\left(\boldsymbol{\mu}, t\right)-\tilde{\mathbf{u}} _{h}\left(\boldsymbol{\mu}, t\right)\right\|_{\ell^2}}{\left\|\mathbf{u}_{h}\left(\boldsymbol{ \mu},t\right)\right\|_{\ell^2}}.
    \label{single}
\end{equation}
The average relative error in the entire test set is defined as
\begin{equation}
    \mathcal{E}_{total}=\frac{1}{N_{\text {test }}} \sum_{i=1}^{N_{test}}\left(\frac{\sum_{k= 1}^{N_{t}}\left\|\mathbf{u}_{h}^{k}\left(\boldsymbol{\mu}_{\text {test }, i}\right)-\tilde{\mathbf{u}}_{h}^{k}\left(\boldsymbol{\mu}_{\text {test }, i}\right)\right\|_{\ell^2}} {\sum_{k=1}^{N_{t}}\left\|\mathbf{u}_{h}^{k}\left(\boldsymbol{\mu}_{\text {test }, i }\right)\right\|_{\ell^2}}\right).
    \label{total}
\end{equation}

\begin{figure}[h]
    \centering
    \includegraphics[width=0.4\textwidth,height=4cm]{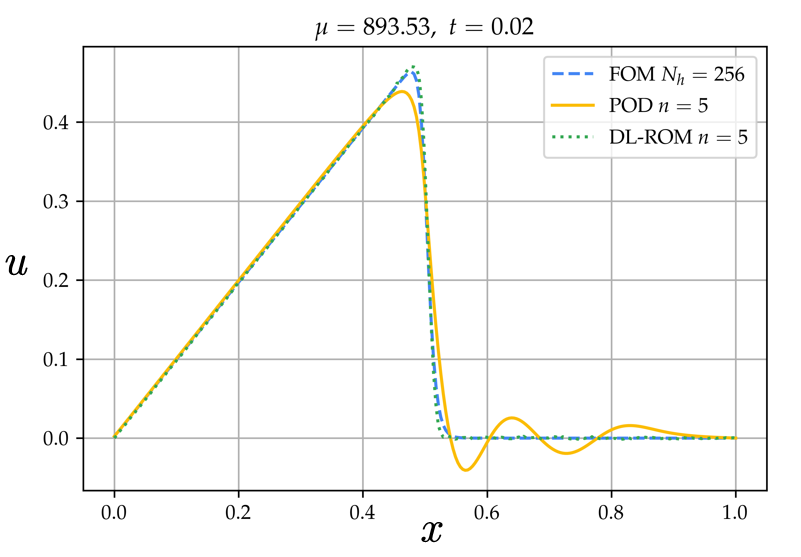}
    \hspace{0.5cm}
    \includegraphics[width=0.4\textwidth,height=4cm]{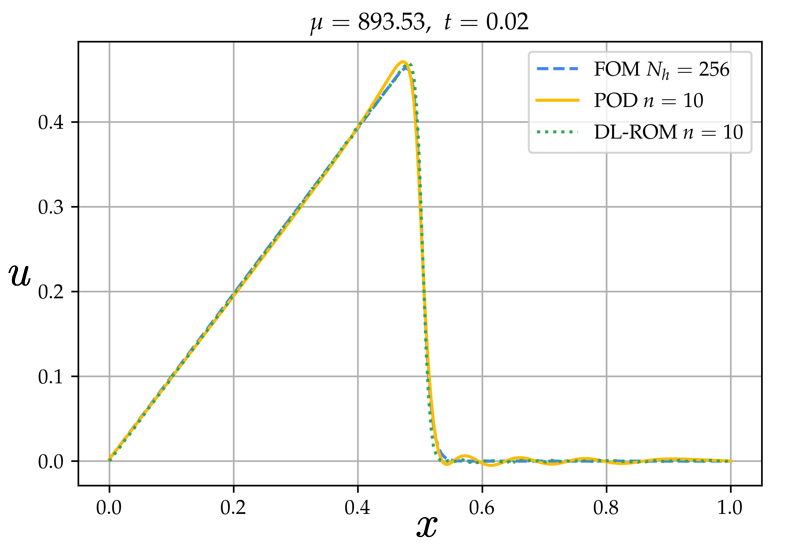}
    \caption{FOM, POD and DL-ROM solutions of parameter $\mu_{\text{test}}= 893.53$ at $t = 0.02$. The dimension of the reduced space: left $n = 5$, right $n=10$.}
    \label{fig:burgers_ad_test_case}
\end{figure}
Figure \ref{fig:burgers_ad_test_case} compares the results of using POD and DL-ROM to solve the equation\,\eqref{burgers} corresponding to parameter $\mu_{\text{test}}= 893.53$ at $t = 0.02$, using the FOM solution as a reference. We can observe that when solving this parameter with the same reduced dimension, DL-ROM has higher accuracy. 
\begin{figure}[h]
    \centering
    \begin{minipage}[t]{0.48\textwidth}
    \centering
    \includegraphics[width=7cm,height=4.5cm]{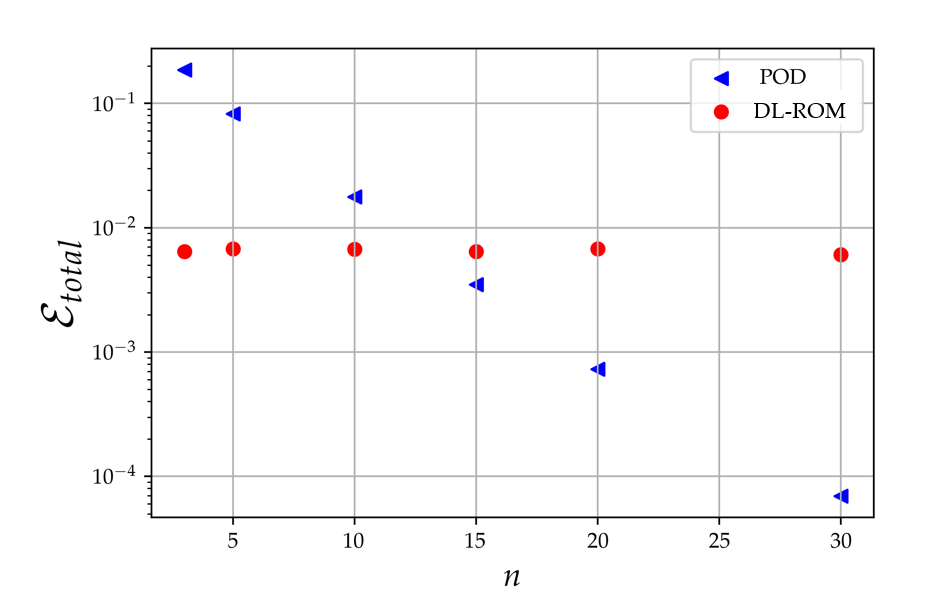}
    \caption{The average relative error of POD and DL-ROM on the test set, where $\mathcal{E}_{total}$ is defined in \eqref{total}.}
    \label{fig:burgers_ad_test_error}
    \end{minipage}
    \begin{minipage}[t]{0.48\textwidth}
        \centering
        \includegraphics[width=8cm,height=5cm]{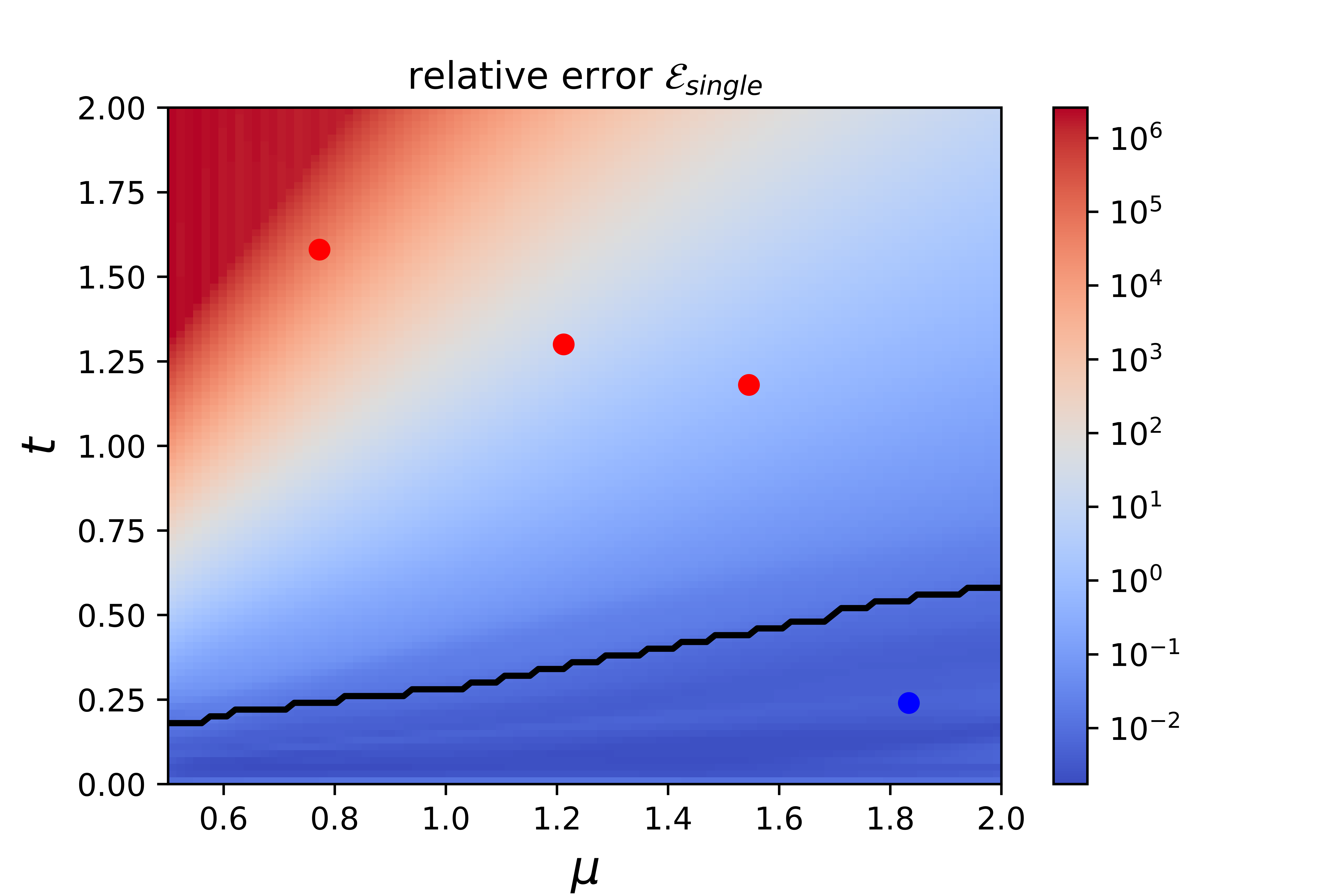}
        \caption{The relative error of DL-ROM solving viscous Burgers equation on the test set.}
        \label{fig:Burgers_diffu_test_error}
    \end{minipage}
\end{figure}
Figure \ref{fig:burgers_ad_test_error} displays the change of the average relative
error of the test set \eqref{total}. We can find that when $n$ is smaller, the
approximation of DL-ROM is better than POD. These results are consistent with
\cite{Fresca2020}. 

\begin{figure}[h]
    \centering
    \includegraphics[width=0.4\textwidth,height=4cm]{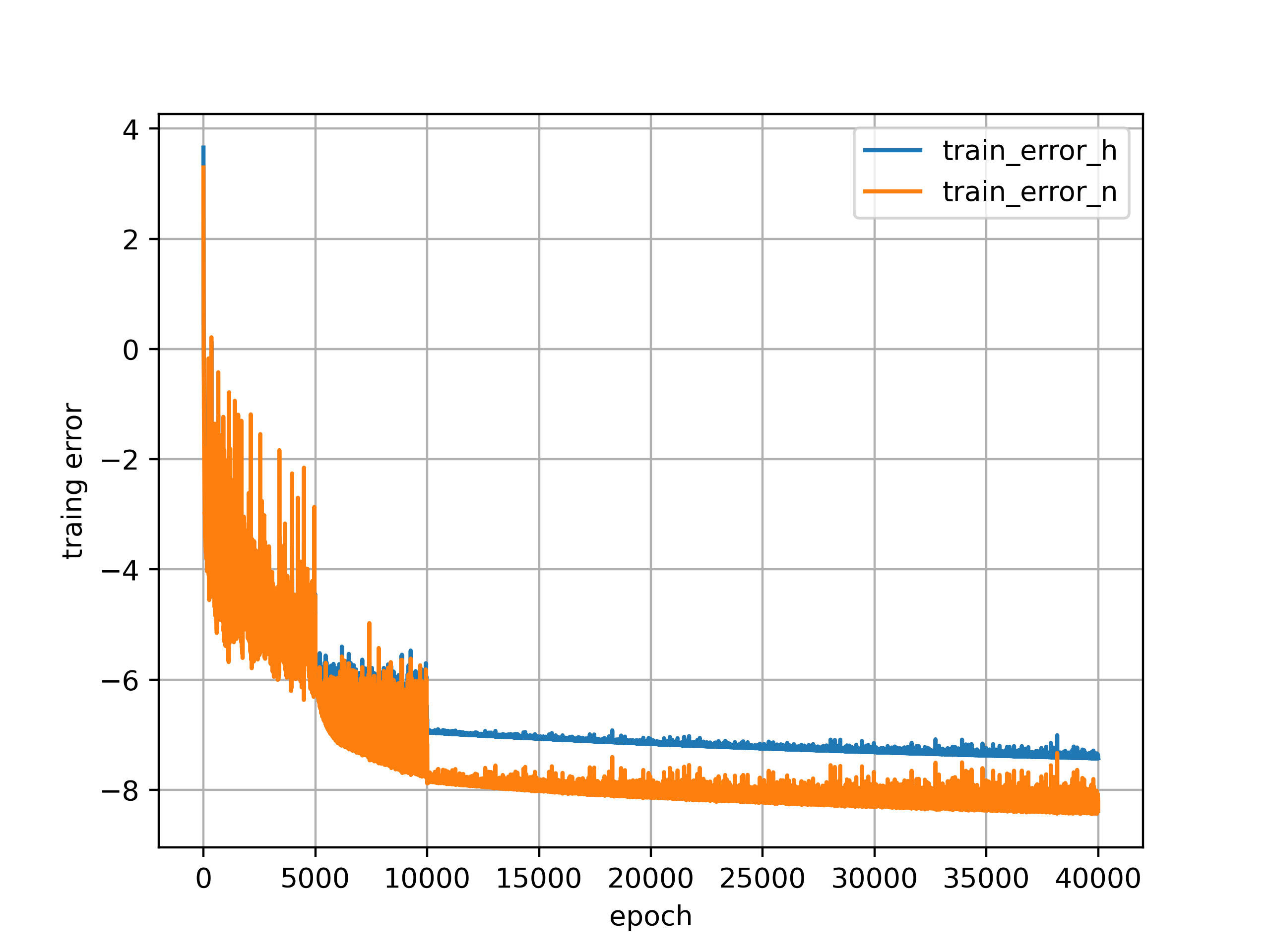}
    \hspace{0.5cm}
    \includegraphics[width=0.4\textwidth,height=4cm]{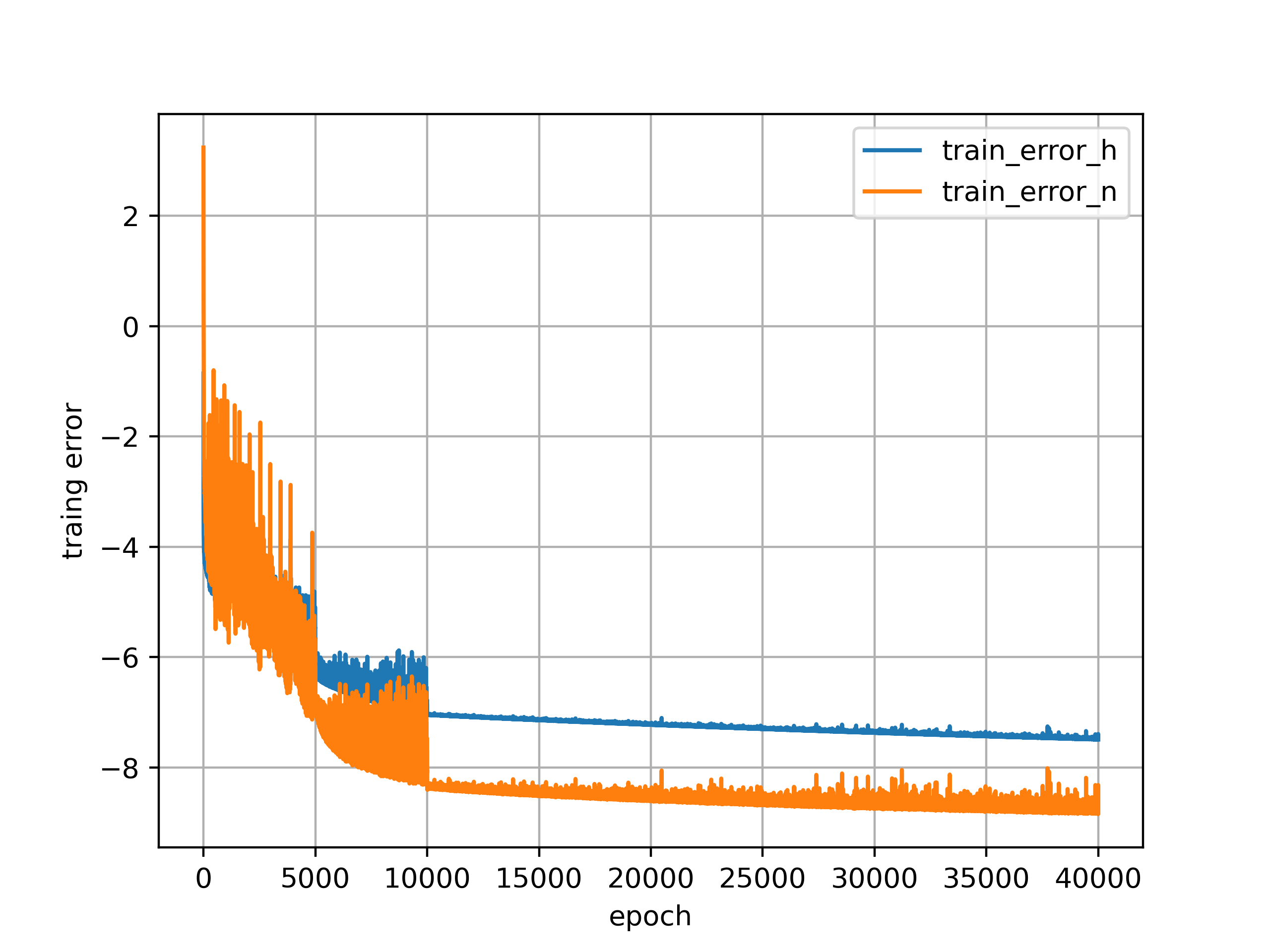}
    \caption{Training error of DL-ROM with 1D (Left) and 2D (Right) convolution layers.
    ``train\_error\_h" and ``train\_error\_n'' refer to the first part of the error and the
    second part of the error in \eqref{eq:316}, respectively. 
	}
    \label{fig:burgers_trainerror}
\end{figure}

It is worth pointing out that in the above experiments, we reshape the $1
\mathrm{D}$ data in $\mathbb{R}^{256}$ to a $2 \mathrm{D}$ data in $\mathbb{R}^{16
\times 16}$ to apply $2 \mathrm{D}$ convolutional neural networks as DL-ROM does.
However, it may be more suitable to directly apply $1 \mathrm{D}$ convolution to $1
\mathrm{D}$ data. We have made a comparison of 1D and 2D convolutional network lays.
Table \ref{tab:Network_1dconv} shows the corresponding network structure parameters.
Figure \ref{fig:burgers_trainerror} gives corresponding results.
The left and right plots of Figure \ref{fig:burgers_trainerror} show the
training error of DL-ROM using 1D convolution and 2D convolution layers.
It can be found that they have similar training error: 0.007444 (Conv1D) and
0.007086 (Conv2D).
Further analysis demonstrates that the 1D convolutional neural networks can
save nearly half network parameters than the 2D networks. 
However, for a fair comparison, we still use the 2D network layers
as Refs. \cite{Lee2020, Fresca2020} did in the following experiments.

\begin{table}[htb]
\centering
\caption{Network architecture of DL-ROM with 1D convolutional and 1D
deconvolutional layers. $N_{b}$ represents batch size. Notations in Conv1d and
ConvTranspose1d are: i: in channels; o: out channels; k: kernel size; s: stride; p:
padding; op: output padding.}
\label{tab:Network_1dconv}

\resizebox{\textwidth}{!}{
\begin{tabular}{|llr|ll|}
    \hline
    \multicolumn{2}{|c}{$\boldsymbol{h}_{\mathrm{enc}}$ (input shape:($N_b$, 256))} &        & \multicolumn{2}{c|}{$\boldsymbol{h}_{\mathrm{dec}}$ (input shape:($N_b$, $n$))} \\
    \hline
    Layer type & Output shape &        & Layer type & Output shape \\
    \hline
    Conv1d(i=1,o=4,k=5,s=2,p=2) & ($N_b$,4, 128) &        & Fully connected & ($N_b$,256) \\
    Conv1d(i=4,o=8,k=5,s=2,p=2) & ($N_b$,8, 64) &          & Fully connected & ($N_b$,256) \\
    Conv1d(i=8,o=16,k=5,s=2,p=2) & ($N_b$,16, 32) &          & Reshape & ($N_b$,16, 16) \\
    Conv1d(i=16,o=16,k=5,s=2,p=2) & ($N_b$,16, 16) &         & ConvTranspose1d(i=16,o=16,k=5,s=2,p=1) & ($N_b$,16,33) \\
    Reshape & ($N_b$,256)  &                  & ConvTranspose1d(i=16,o=8,k=5,s=2,p=2) & ($N_b$,8,65) \\
    Fully connected & ($N_b$,256) &          & ConvTranspose1d(i=8,o=4,k=5,s=2,p=2) & ($N_b$,4,129) \\
    Fully connected & ($N_b$,$n$)  &         & ConvTranspose1d(i=4,o=1,k=5,s=2,p=3,op=1) & ($N_b$,1,256) \\
    \hline
\end{tabular}
}
\end{table}

Next we change the parameter range to $\mathcal{P}=[0.5,2]$, and its corresponding viscosity coefficient $1/\mu$ belongs to $[0.5,2]$. We also sample $N_{\text {train}}=20$ parameters in $\mathcal{P}$ according to the uniform distribution to generate the training set ($N_h=256, N_t=100$). The test set is obtained by sampling $N_{\text {test}}=100$ parameters equidistantly in the parameter range.

Figure \ref{fig:Burgers_diffu_test_error} shows the relative error of the DL-ROM's prediction results for each parameter in the test set. The relative error above the black line is greater than $10^{-2}$, and below is less than $10^{-2}$. 
\begin{figure}[h]
    \centering
    \subfigure{\includegraphics[width=0.4\textwidth,height=4.0cm]{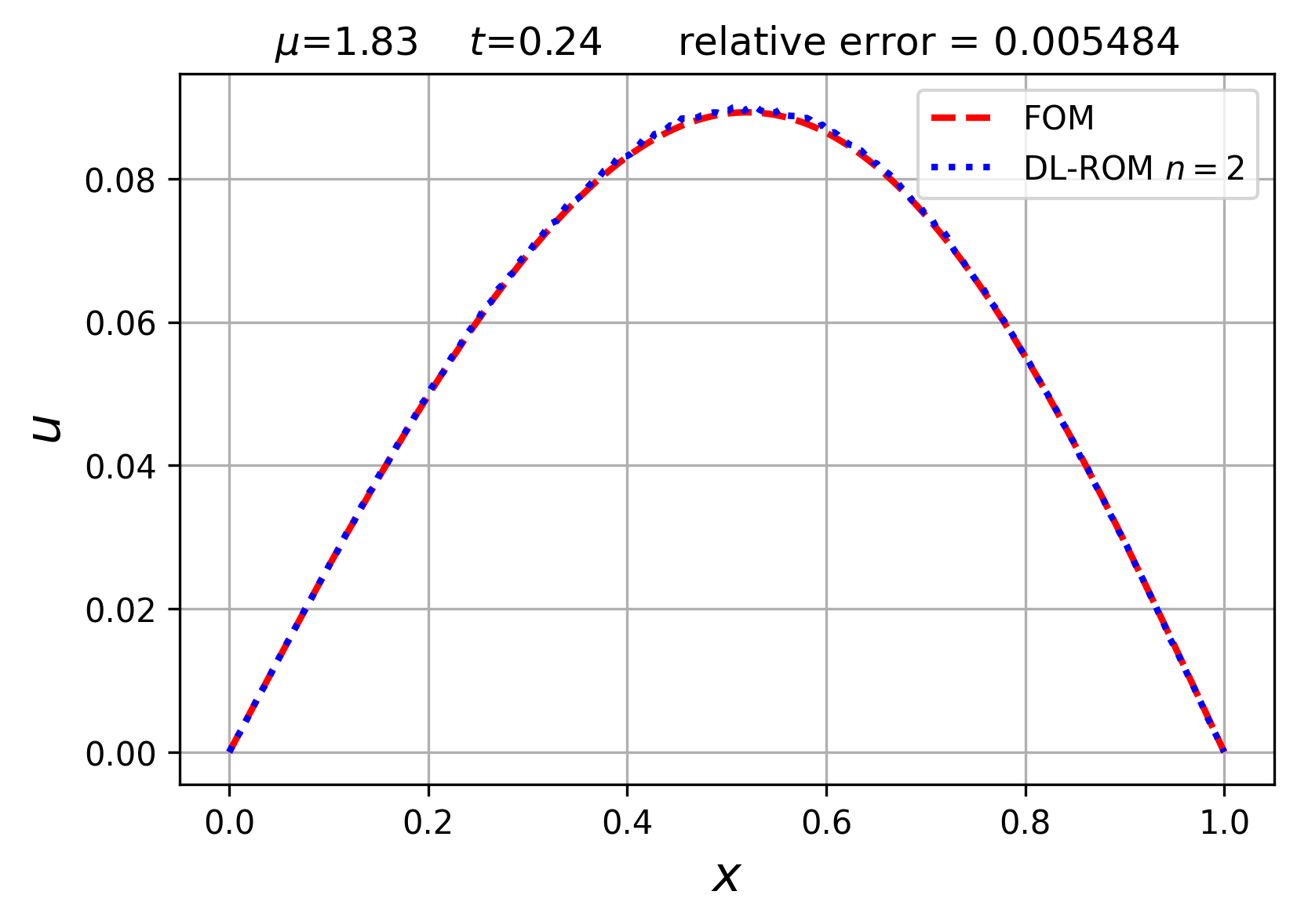}}
    \subfigure{\includegraphics[width=0.4\textwidth,height=4.0cm]{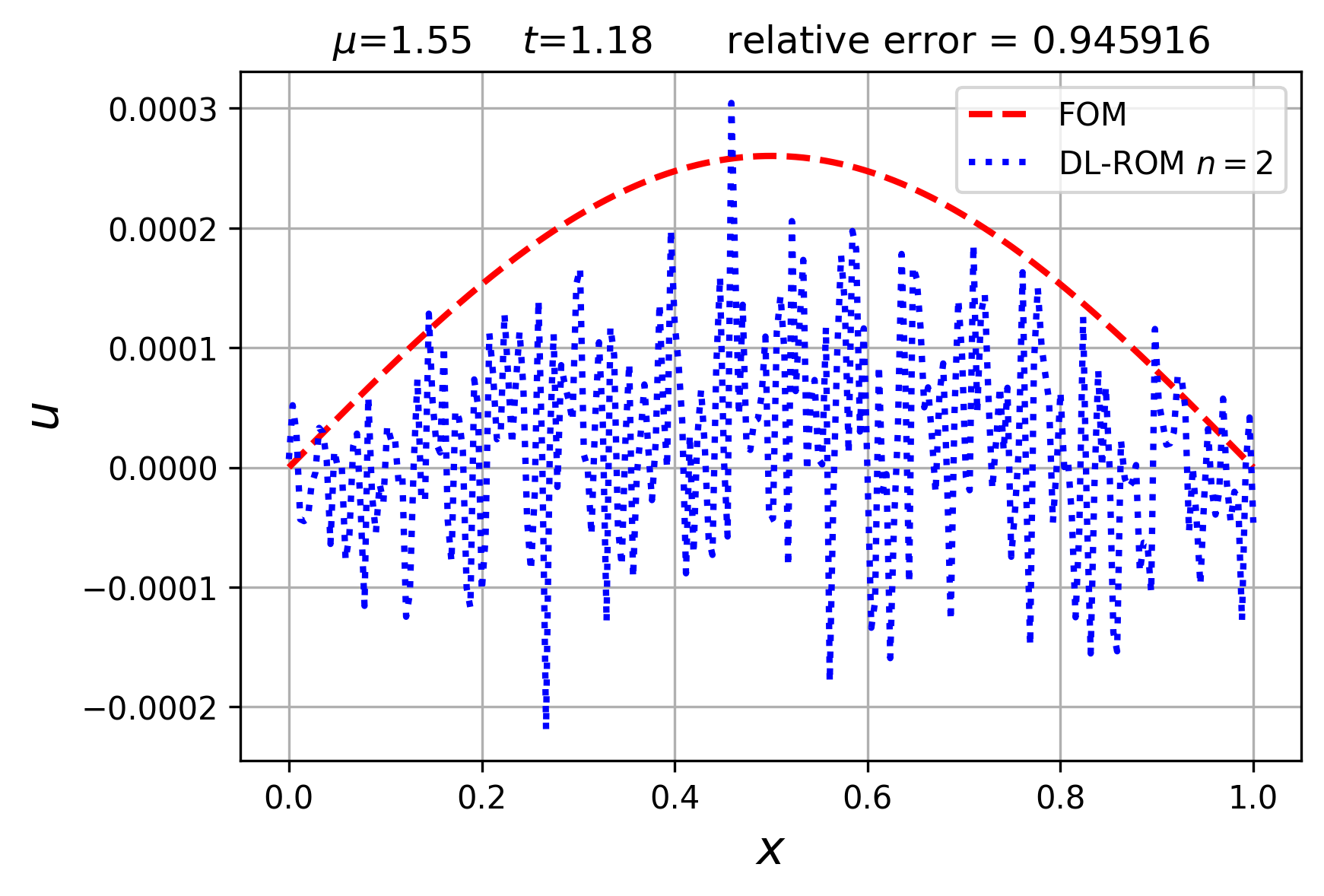}}
    \\ 
    \centering
    \subfigure{\includegraphics[width=0.4\textwidth,height=4.0cm]{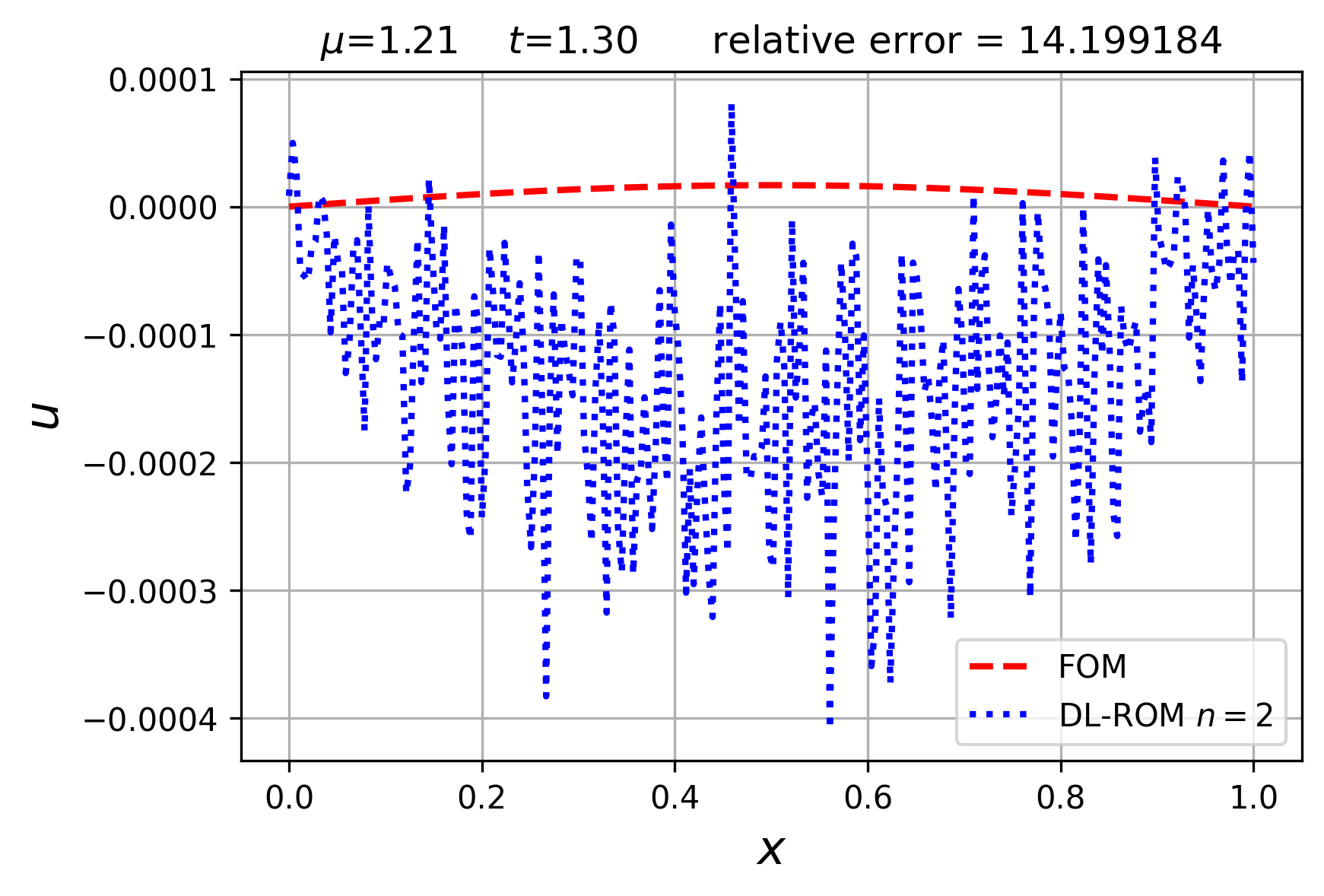}}
    \subfigure{\includegraphics[width=0.4\textwidth,height=4.0cm]{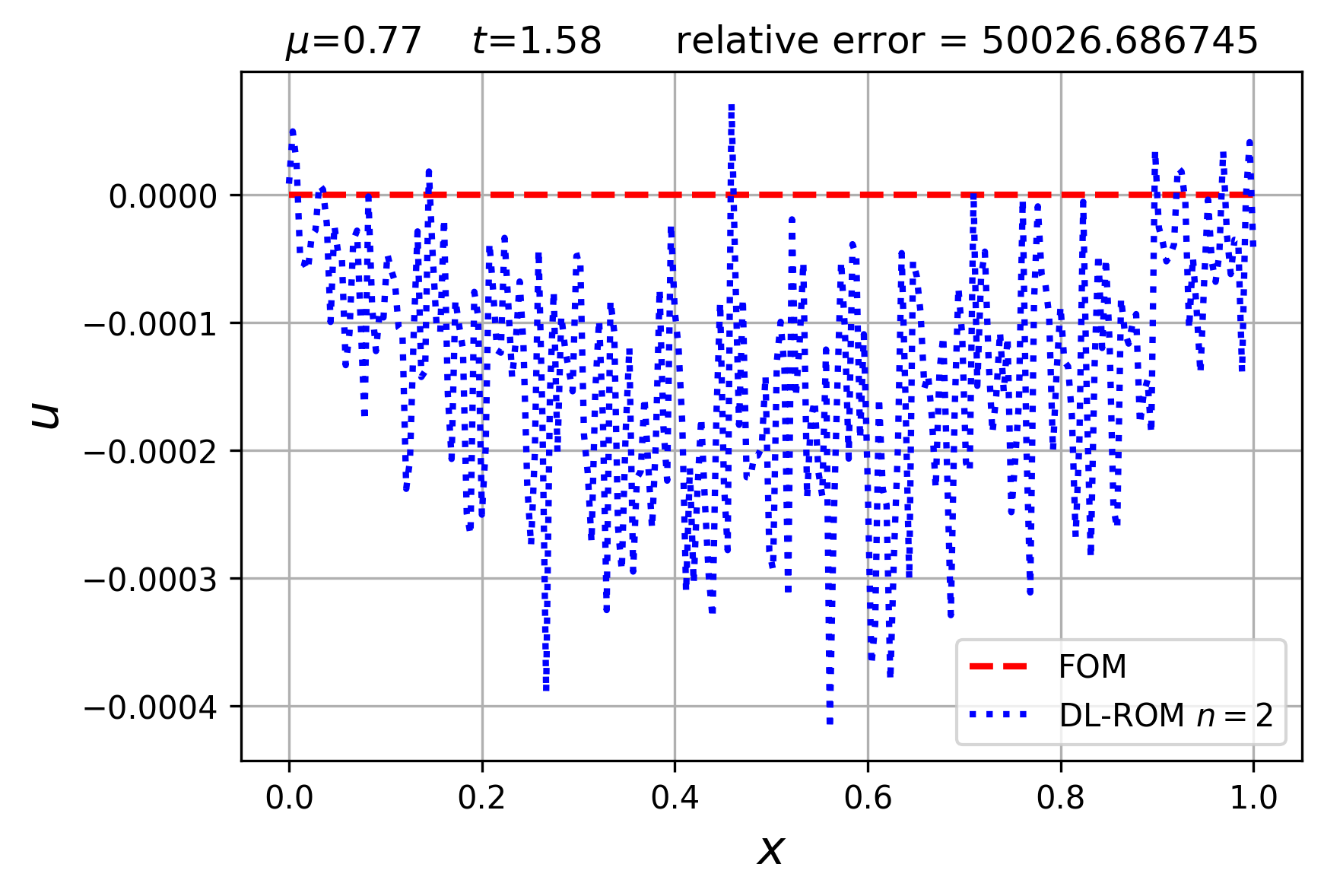}}
    \captionsetup{justification=centering}
    \caption{The results of viscous Burgers equation corresponding to the marked parameters in Figure \ref{fig:Burgers_diffu_test_error}.}
    \label{fig:burgers_diffu_test_case}
\end{figure}
Figure \ref{fig:burgers_diffu_test_case} shows DL-ROM's inference results of the parameters corresponding to the four marked points in Figure \ref{fig:Burgers_diffu_test_error}. This implies that DL-ROM's generalization ability becomes weakened.

The reason is attributed to the fact that when the parameter range is set to $\mathcal{P}=[0.5, 2]$, the diffusion term in the equation \eqref{burgers} plays a role, causing its solution to change drastically along the time direction. We define the following $\gamma$ to quantify the severity of the solution change

\begin{equation}
    \gamma=\frac{\max _{\mu \in \mathcal{P} \atop t \in[0,
	T]}\left\|\boldsymbol{u}_{\boldsymbol{h}}(x, t,
	\boldsymbol{\mu})\right\|_{\ell^{\infty}}}{\min _{\mu \in \mathcal{P} \atop t
	\in[0, T]}\left\|\boldsymbol{u}_{h}(x, t,
	\boldsymbol{\mu})\right\|_{\ell^{\infty}}}.
\end{equation}
The larger the $\gamma$, the more drastic the solution changes, and vice versa. For viscous Burgers equation \eqref{burgers}, when $\mathcal{P}=[100,1000] $, $\log_{10}\gamma <2$, the solution does not change much and DL-ROM has good generalization ability. On the contrary, when $\mathcal{P}=[0.5, 2]$, $\log_{10}\gamma> 2$, the solution changes drastically, DL-ROM is inclined to fit the solution whose norm is relatively large but ignore solutions with relatively small norms. This observation inspires us to classify the original dataset and establish a classification network based on the classified data.

\section{Multiclass classification-based ROM}\label{sec05}

Inspired by the observation above, in this section, we design the MC-ROM and analyze its computational complexity of online computation.

\subsection{Network architecture}

\begin{figure}[h]
    \centering
    \includegraphics[width=0.9\textwidth]{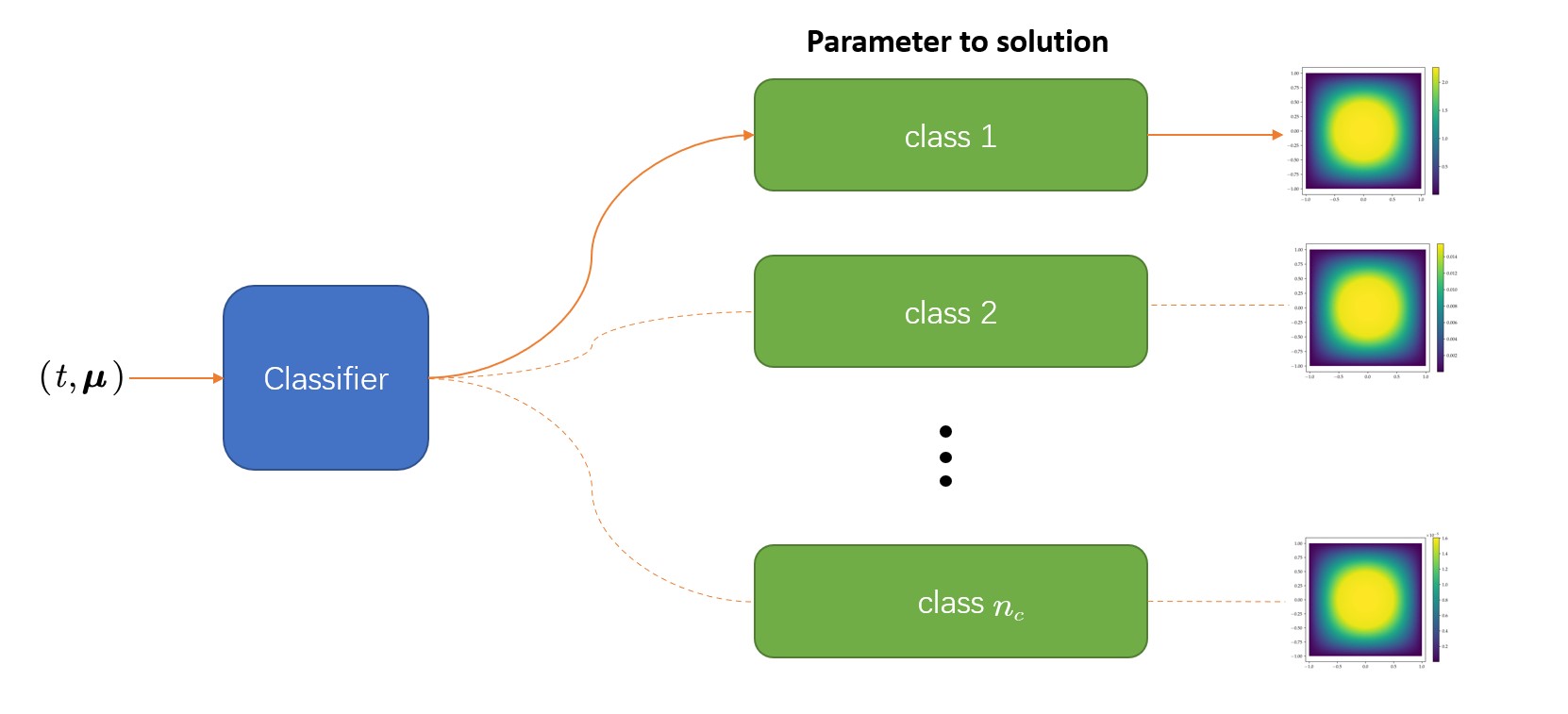}
    \captionsetup{justification=centering}
    \caption{MC-ROM network architecture.}
    \label{fig:MC_ROM}
\end{figure}

Figure \ref{fig:MC_ROM} gives the network architecture of MC-ROM. The blue block on the left plot is the classifier, and its purpose is to classify the parameters according to the order of magnitude of the numerical solution. Its input is $(t,\boldsymbol{\mu})$, and the output is the label to classify the data. In practice, we use the order of magnitude of $\ell^\infty$ norm of numerical solution as classification criteria to assign each parameter a label. We then use parameter-label pairs as the training dataset for the classifier.

\begin{figure}[h]
     \centering
     \includegraphics[width=0.5\textwidth]{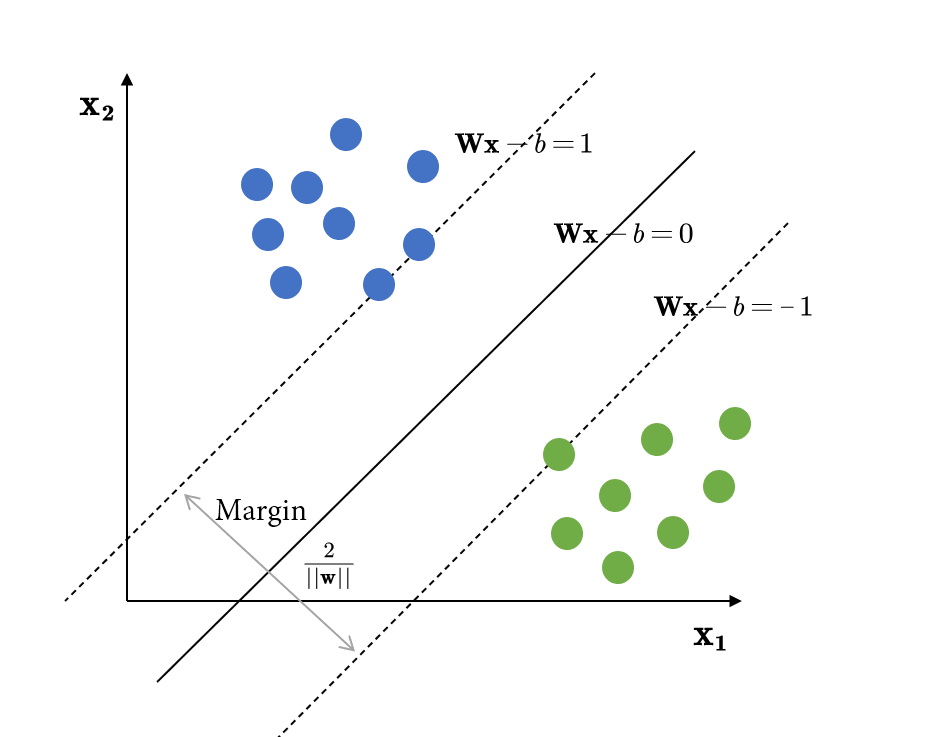}
     \captionsetup{justification=centering}
     \caption{SVM is to maximum the margin.}
     \label{fig:SVM}
\end{figure}

There are many classification algorithms, such as neural networks, decision trees,
k-nearest neighbors, naive Bayes, etc. The choice of classification algorithm depends
on the problem. Here, we choose SVM to train the classifier due to the small training
data. For our numerical results, SVM has better performance than fully connected
neural networks. As shown in Figure \ref{fig:SVM}, SVM separates different types of
data by constructing a hyperplane or a group of hyperplanes in a high-dimensional or
infinite-dimensional space. The distance between different classes is called margin.
Generally, the larger the margin, the smaller the generalization error of the
classifier, so the hyperplane should be constructed so that the distance from the
nearest training data point of any class is large enough. SVM applies the kernel
trick to make the data set linearly separable and uses Higen loss to achieve these
hyperplanes to meet the above properties. See Refs. \cite{cortes1995support, weston1999support} for a detailed description. We implement SVM through \lstinline{sklearn}\cite{scikit-learn}.

The green blocks on the right plot are different subnets. Each subnet is used to fit the parameter to solution mapping in a certain class. The subnet architectures can be the same or different, depending on the problem. Here, all of our subnet models are DL-ROM due to its high efficiency of online computing and its applicability to convection problems. The offline training process of different subnets are similar learning tasks, and we can apply transfer learning \cite{caruana1997multitask} to speed up training. In transfer learning, when the training of a subnet is completed, its parameters can assign to the next subnet as the initial parameters for training. The numerical experiment in Section \ref{sec61} will show the advantages of transfer learning.
The offline training and online testing of MC-ROM are summarized in Algorithm \ref{algorithm:MC_ROM_offline} and Algorithm \ref{algorithm:MC_ROM_online}.

\begin{algorithm}[htb]
    \caption{MC-ROM offline traing}\label{algorithm:MC_ROM_offline}
    \KwData{Parameter matrix $P \in \mathbb{R}^{\left(n_{\mu}+1\right) \times N_{s}},$ snapshot matrix $S \in \mathbb{R}^{N_{h} \times N_{s}}$}
    \KwIn{The number of training epochs $N_{\text {epochs}},$ batch size $N_{b}$, learning rate $\eta$, hyperparameter $\alpha, \beta $} 
    \KwOut{A classifier $\mathcal{C}$ and trained models $\boldsymbol{h}_{\mathrm{enc}}^i,\boldsymbol{h}_{\mathrm{dec}}^i,\Psi^i,$ $i=1,\cdots,n_{\text{c}}$} 
    Classify $P$ and $S$ into $[P_1,\dots,P_{n_c}]$ and $[S_1,\dots,S_{n_c}] $ and get label matrix $L \in \mathbb{R}^{n_{\text{c}}\times N_{s}}$
 
    Use SVM to train a classifier $\mathcal{C}$ with $P$ and $L$

    \For{$i=1:n_{\text{c}}$}{
        Substitute the data $P_i, S_i$ into Algorithm \ref{algorithm:DL_ROM_offline} and use DL-ROM input parameters to train to get $\boldsymbol{h}_{\mathrm{enc}}^i,\boldsymbol{h}_{\mathrm{dec}}^i,\Psi^i$ 
    }
    \Return{$\mathcal{C}$ and $\boldsymbol{h}_{\mathrm{enc}}^i,\boldsymbol{h}_{\mathrm{dec}}^i,\Psi^i,\,i=1,\cdots,n_{\text{c}}$}
\end{algorithm}

\begin{algorithm}[htb]
    \caption{MC-ROM online testing}\label{algorithm:MC_ROM_online}
    \KwData{Parameter matrix $P \in \mathbb{R}^{\left(n_{\mu}+1\right) \times N_{s}}$}
    \KwIn{Algorithm \ref{algorithm:MC_ROM_offline} trained model $\mathcal{C}$ and $\Psi^i,\boldsymbol{h}_{\mathrm{dec}}^i;\,\,i=1,\cdots,n_{\text{c}}$}
    \KwOut{The predicted solution matrix $\widetilde{S}$ corresponding to the parameter matrix $P$}
    Substitute the parameters into the classifier $\mathcal{C}$ to determine which class it belongs to

    Substitute the parameters into the corresponding trained model $\boldsymbol{h}_{\mathrm{dec}}^i\circ \Psi^i$ to get $\widetilde{S}$
\end{algorithm}

\subsection{Computational complexity}\label{sec52}

This subsection analyzes the online computational complexity of MC-ROM and POD. We make a convention that
\begin{equation}
     \mathbb{A}_{h}(\boldsymbol{\mu}) \mathbf{u}_{h}(\boldsymbol{\mu})=\mathbf{f}_{h}(\boldsymbol{\mu }),
     \label{eq:fom}
\end{equation}
represents the FOM system, where $\mathbb{A}_{h}(\boldsymbol{\mu}) \in \mathbb{R}^{N_{h} \times N_{h}}$, and
\begin{equation}
     \mathbb{A}_{n}(\boldsymbol{\mu}) \mathbf{u}_{n}(\boldsymbol{\mu})=\mathbf{f}_{n}(\boldsymbol{\mu }),
     \label{eq:rom}
\end{equation}
represents the reduced system, $\mathbb{A}_{n}(\boldsymbol{\mu}) \in \mathbb{R}^{n\times n}$, $n$ may be much smaller than $N_h$.

The online calculation of MC-ROM has two steps except classification:
\begin{enumerate}
    \item Solve $\mathbf{u}_{n}$. Use a fully connected network to approximate $\mathbf{u}_{n}(\boldsymbol{\mu})$, the computional amount of each layer is $N_{\text{neu}}^2$. If there are $l$ hidden layers, the computional amount is $l\times N_{\text{neu}}^2$ except the input and output layers;
    \item Lift $\mathbf{u}_n$ to high-dimensional manifold through the decoder. The computional amount of a deconvolution layer is $O\left(N_h \cdot k^{2} \cdot C_{ \text {in }} \cdot C_{\text {out }}\right)$, where $k^2$ is the size of the convolution kernel. $C_{\text {in }}$ and $C_{\text {out}}$ represent the number of channels of the previous layer and the current layer.
\end{enumerate}
In actual calculations, $N_{\text{neu}},\,l,\,k,\,C_{\text {in}},\,C_{\text {out}}$ are much smaller than $N_h$, then the online computational complexity of MC-ROM is $O(N_h)$. 

The online stage of POD can has the following four steps:
\begin{enumerate}
    \item Given a $\boldsymbol{\mu}$, discrete its corresponding PDE to get \eqref{eq:fom};
    \item Project the linear system \eqref{eq:fom} to a low-dimensional manifold and get \eqref{eq:rom}, where
    \begin{equation}
        \mathbb{A}_{n}=V^{T} \mathbb{A}_{h} V, \quad \mathbf{f}_{n}=V^{T} \mathbf{f}_{h}.
    \end{equation}
    and $V\in \mathbb{R}^{N_h \times n}$ is obtained in POD offline stage;
    \item Solve \eqref{eq:rom};
    \item Lift the low-dimensional solution $\mathbf{u}_{n}(\boldsymbol{\mu})$ to high-dimensional manifold $V\mathbf{u}_{n}(\boldsymbol{\mu}) $.
\end{enumerate}
The computational complexity corresponding to each step above is as follows:
\begin{enumerate}
    \item Use spatial discretization methods such as FEM, FDM to obtain \eqref{eq:fom}, computational complexity depends on the discretization methods used;
    \item Projecte it to a low-dimensional manifold requires matrix multiplication, hence the computational complexity is $O\left({N_h}^2\right)$;
    \item Equation \eqref{eq:rom} is dense, the computational complexity of directly solving it is $O\left(n^{3}\right)$;
    \item Lift operation is $V\mathbf{u}_n(\boldsymbol{\mu})$, the operations required is $N_h\times n$.
\end{enumerate}
Evidently, the online computational complexity of the POD is much large than the MC-ROM.

Sometimes, the first two parts can be put offline if
$\mathbb{A}_{h}(\boldsymbol{\mu})$ and $\mathbf{f}_{h}(\boldsymbol{\mu})$ in
\eqref{eq:fom} are \textit{parameter-separable}, 
\begin{equation}
    \mathbb{A}_{h}(\boldsymbol{\mu})=\sum_{q=1}^{Q_{a}} \theta_{a}^{q}(\boldsymbol{\mu}) \mathbb{A}_{h}^{q}, \quad \mathbf{f}_{h}(\boldsymbol{\mu})=\sum_{q=1}^{Q_{f}} \theta_{ f}^{q}(\boldsymbol{\mu}) \mathbf{f}_{h}^{q},
    \label{fom_affine}
\end{equation}
where $\theta_{a}^{q}: \mathcal{P} \rightarrow \mathbb{R}, \,q=1, \ldots, Q_{a}$ and
$\theta_{f}^{q}: \mathcal{P} \rightarrow \mathbb{R}, \,q=1, \ldots, Q_{f}$.
Then we can apply the Galerkin method \eqref{ROM} to \eqref{fom_affine} and obtain
\begin{equation}
    \mathbb{A}_{n}(\boldsymbol{\mu})=\sum_{q=1}^{Q_{a}} \theta_{a}^{q}(\boldsymbol{\mu}) \mathbb{A}_{n}^{q}, \quad \mathbf{f}_{n}(\boldsymbol{\mu})=\sum_{q=1}^{Q_{f}} \theta_{ f}^{q}(\boldsymbol{\mu}) \mathbf{f}_{n}^{q},
    \label{rom_affine}
\end{equation}
where $\mathbb{A}_{n}^{q}=V^{T} \mathbb{A}_{h}^{q} V, \,\, \mathbf{f}_{n}^{q }=V^{T} \mathbf{f}_{h}^{q}$ can be pre-calculated. In this case, the first two steps of POD online calculation are replaced by directly assembling equation \eqref{rom_affine}. The operations required to generate \eqref{rom_affine} is
$O\left(Q_{a} n^{2}+Q_{f} n\right)$. If $n$ and $Q_{a},\,Q_f$ are small enough, then the online computational complexity of POD is $O(N_h)$. 

However, even for the parameter-separable problems, MC-ROM still has the following
advantages over POD. Firstly, we noticed that the computational complexity of the
third step of POD is $O\left(n^{3}\right)$. Therefore, when the problem to be solved
requires a high-dimensional reduced space, the online calculation of POD will become
expensive. Secondly, using neural networks as the surrogate model to establish the mapping between parameters and solutions has fundamentally changed the calculation process of many scenarios in scientific computing. Compared with POD, MC-ROM is a non-intrusive ROM. The calculation based on the neural network makes it easier to use computing resources such as GPU, making its online calculation much faster than POD. Section \ref{sec06} shows the related numerical experiments.

\section{Numerical Experiments}\label{sec06}

This section applies MC-ROM to solve the $1 \mathrm{D}$ viscous Burgers equation and 2D parabolic equation with discontinuous diffusion coefficients. We also compare the approximation accuracy and computational time of solving these equations with DL-ROM and with POD. All the experiments in this section are done on a workstation equipped with two Intel(R) Xeon(R) Silver 4214 CPU @ 2.20GHz, 128GB RAM, and Nvidia Tesla V100-PCIe-32GB GPU.

\subsection{$1 \mathrm{D}$ viscous Burgers equation}\label{sec61}

In Section \ref{sec04}, we have used DL-ROM to solve the viscous Burgers equation\,\eqref{burgers} with parameter interval $\mathcal{P}=[0.5, 2]$, but its generalization ability is very poor. Here, we use MC-ROM to solve it. Firstly, we classify the original training dataset according to the FOM solution vectors. Concretely, we use two consecutive orders of magnitude of $\ell^\infty$ norm to label the parameters.
We divide the original dataset into $n_{\text{c}}$ classes, and the labels corresponding to all parameters are assembled into a label matrix $L \in \mathbb{R}^{n_{\text{c}}\times N_{s}}$. Table \ref{tab:Burgers_clf} shows the classification results of the training set.
\begin{table}[htb]
	\centering
    \caption{Classification results of the training set whose data are FOM solutions of $1 \mathrm{D}$ viscous Burgers equation.}
     \resizebox{\textwidth}{!}{
	\begin{tabular}{|l|c|c|c|c|c|c|}
	\hline
	Range of $\left\lVert\mathbf{u}_{h} \right\rVert_{\ell^{\infty}}$  &   $\left\lVert\mathbf{u}_{h}\right\rVert_{\ell^{\infty}}\geq 10^{-2}$  &   $[10^{-4}, 10^{-2})$     &  $[10^{-6}, 10^{-4})$       &  $[10^{-8}, 10^{-6})$       &   $[10^{-10}, 10^{-8})$      & $\left\lVert\mathbf{u}_{h}\right\rVert_{\ell^{\infty}}\leq 10^{-10}$ \\
	\hline
	Label  & 1      & 2      & 3      & 4      & 5      & 6 \\
	\hline
	Size     & 435    & 573    & 477    & 224    & 244    & 67 \\
	\hline
	\end{tabular}}
	\label{tab:Burgers_clf}
\end{table}

Next we apply the dataset $\left\{p_i, l_i\right\}_{i=1}^{N_s}$ to train the SVM classifier, where $p_i \, ,l_i$ are the $i$-th column of parameter matrix $P \in \mathbb{R}^{\left(n_{\mu}+1\right) \times N_{s}}$ and label matrix $L \in \mathbb{R}^{n_{\text{c}}\times N_{s}}$, respectively.

\begin{figure}[h]
    \centering
    \includegraphics[width=0.5\textwidth]{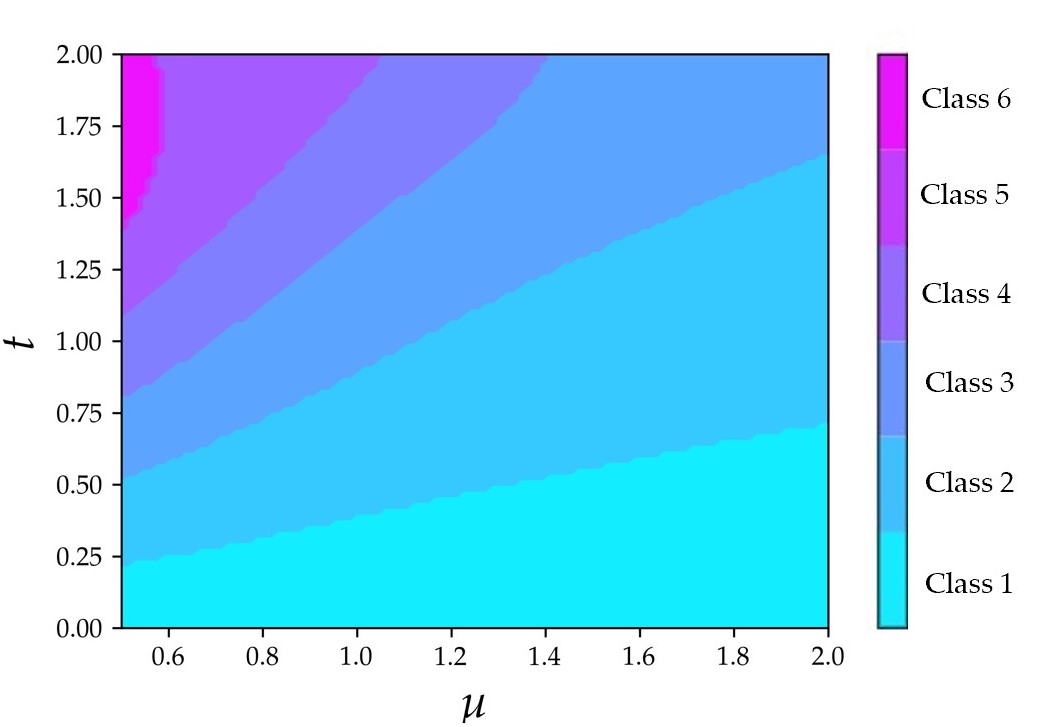}
    \captionsetup{justification=centering}
    \caption{The classification result of the SVM classifier.}
    \label{fig:Burgers_clf}
\end{figure}

We apply the radial basis function kernel in the SVM
\begin{equation}
    K\left(\mathbf{x}, \mathbf{x}^{\prime}\right)=\exp \left(-\gamma\left\|\mathbf{x}-\mathbf{x}^{\prime}\right\|^{2}\right),
\end{equation}
where $\gamma = \frac{1}{(n_{\mu}+1)Var(P)}$, and $Var(P)$ refers to the variance of parameters.
After training, we verify the performance of the classifier on the test set, and the result is that it has an accuracy of $97.17\%$, as shown in Figure \ref{fig:Burgers_clf}. 
And we can see that the black line in the Figure \ref{fig:Burgers_diffu_test_error} is very similar to the line between the first and second class in the Figure \ref{fig:Burgers_clf}. These results demonstrate DL-ROM only approximates well for the first data class. In order to improve DL-ROM so that it has generalization ability for the entire parameter space, it is necessary to train with different subnets according to the classification result shown in Figure \ref{fig:Burgers_clf}.
\begin{figure}[h]
    \centering
    \includegraphics[width=0.5\textwidth]{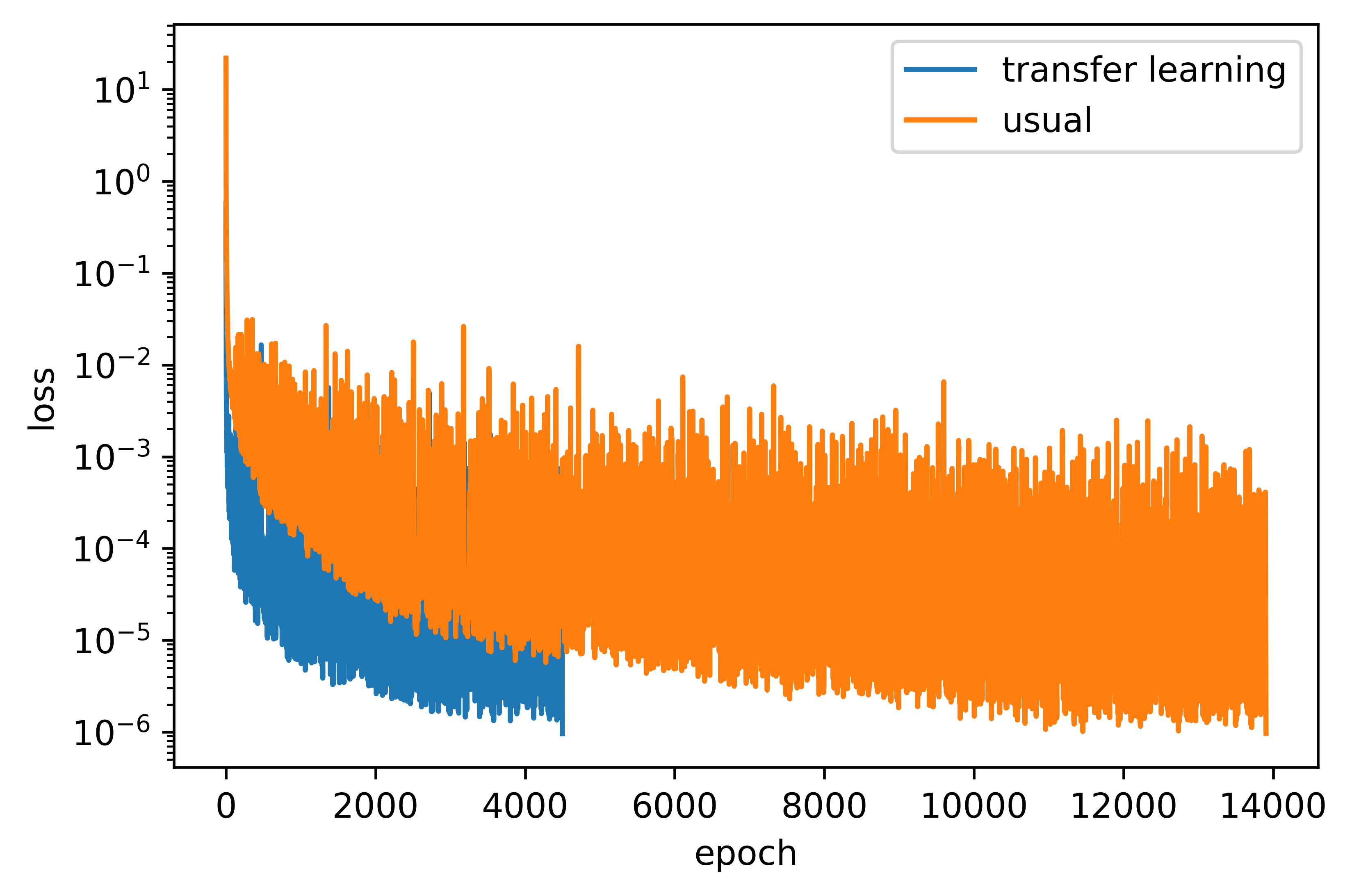}
    \captionsetup{justification=centering}
    \caption{Comparison of the training process with and without transfer learning.}
    \label{fig:transfer}
  \end{figure}

Finally, we use the DL-ROMs as subnets in the MC-ROM which also have the same subnet
architecture as DL-ROM in Ref.\cite{Fresca2020} except for different learning rates which is adjusted during training.  We apply Algorithm \ref{algorithm:DL_ROM_offline} to train each subnet. As mentioned in Section \ref{sec05}, we use transfer learning techniques to accelerate the offline training process. Figure \ref{fig:transfer} compares the impact of using transfer learning and randomly initialized network parameters during offline training. It can be seen that the former is more efficient.

\begin{figure}[h]
  \centering
  \includegraphics[width=0.32\textwidth]{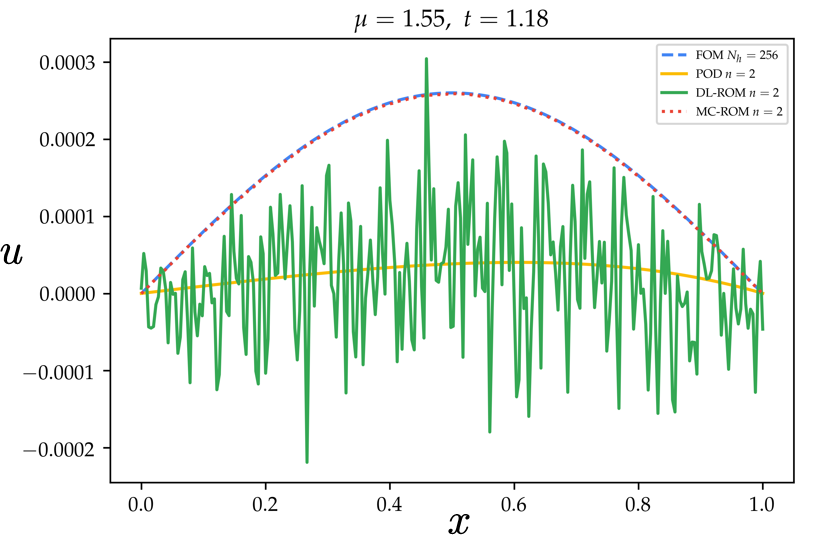}
  \includegraphics[width=0.32\textwidth]{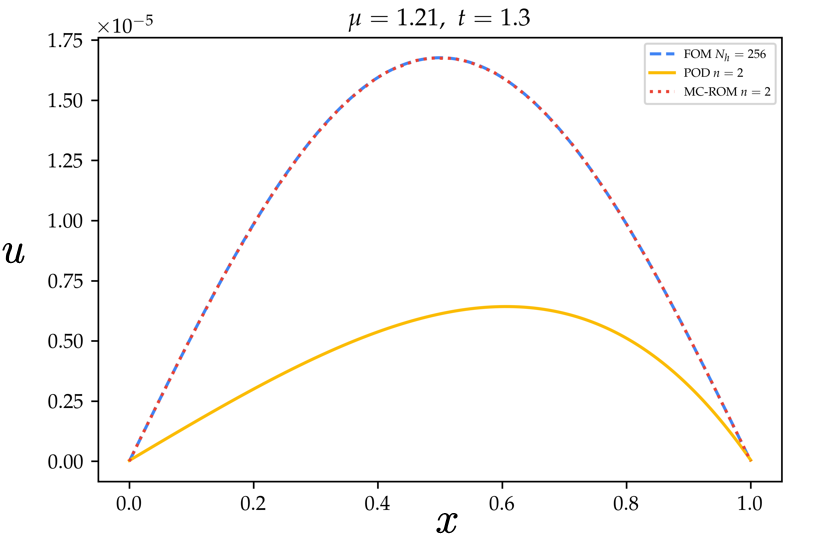}
  \includegraphics[width=0.32\textwidth]{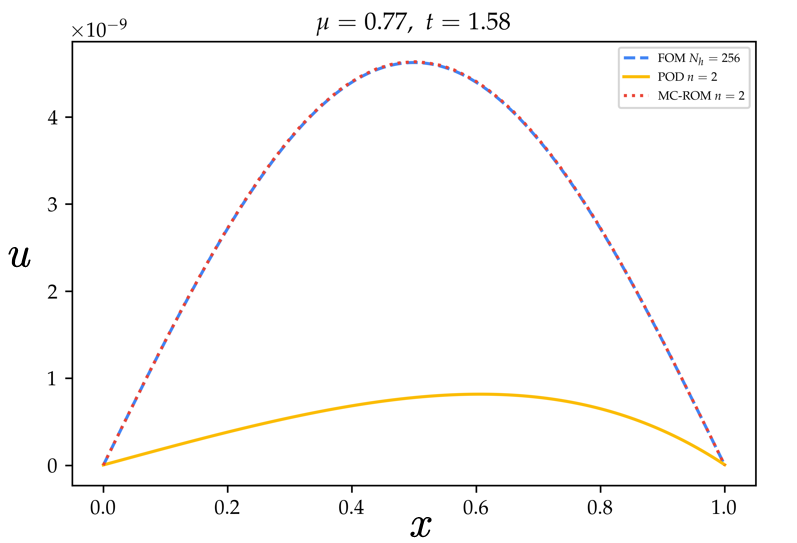}
  \caption{Results of Solving equation \eqref{burgers} corresponding to the red marked parameters in Figure \ref{fig:Burgers_diffu_test_error} by  using POD, DL-ROM, and MC-ROM. $n$ represents the dimension of the reduced space.}
  \label{fig:burgers_diffu_clf_test_case}
\end{figure}

Figure \ref{fig:burgers_diffu_clf_test_case} compares the test results of MC-ROM, DL-ROM, and POD. The numerical results of DL-ROM in the last two plots are not shown since its relative error is too large. Compared to DL-ROM and POD, MC-ROM has good generalization ability for the parameters corresponding to solutions with different orders of magnitude.

\begin{figure}[h]
    \centering
    \subfigure[test error]{
        \label{subfig:a}
        \includegraphics[width=0.45\textwidth]{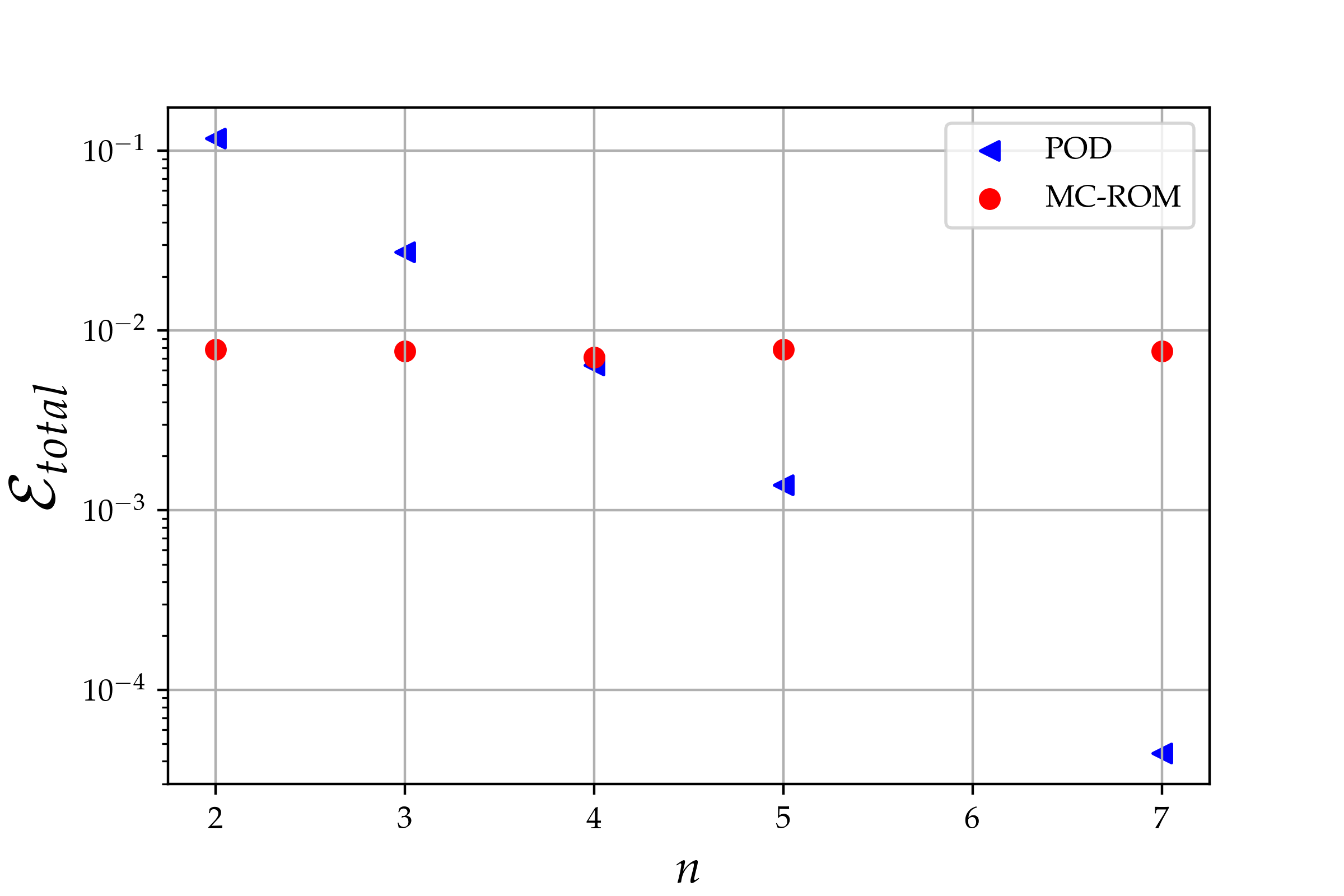}
    }
    \quad
    \subfigure[inference time]{
            \label{subfig:b}
            \includegraphics[width=0.45\textwidth]{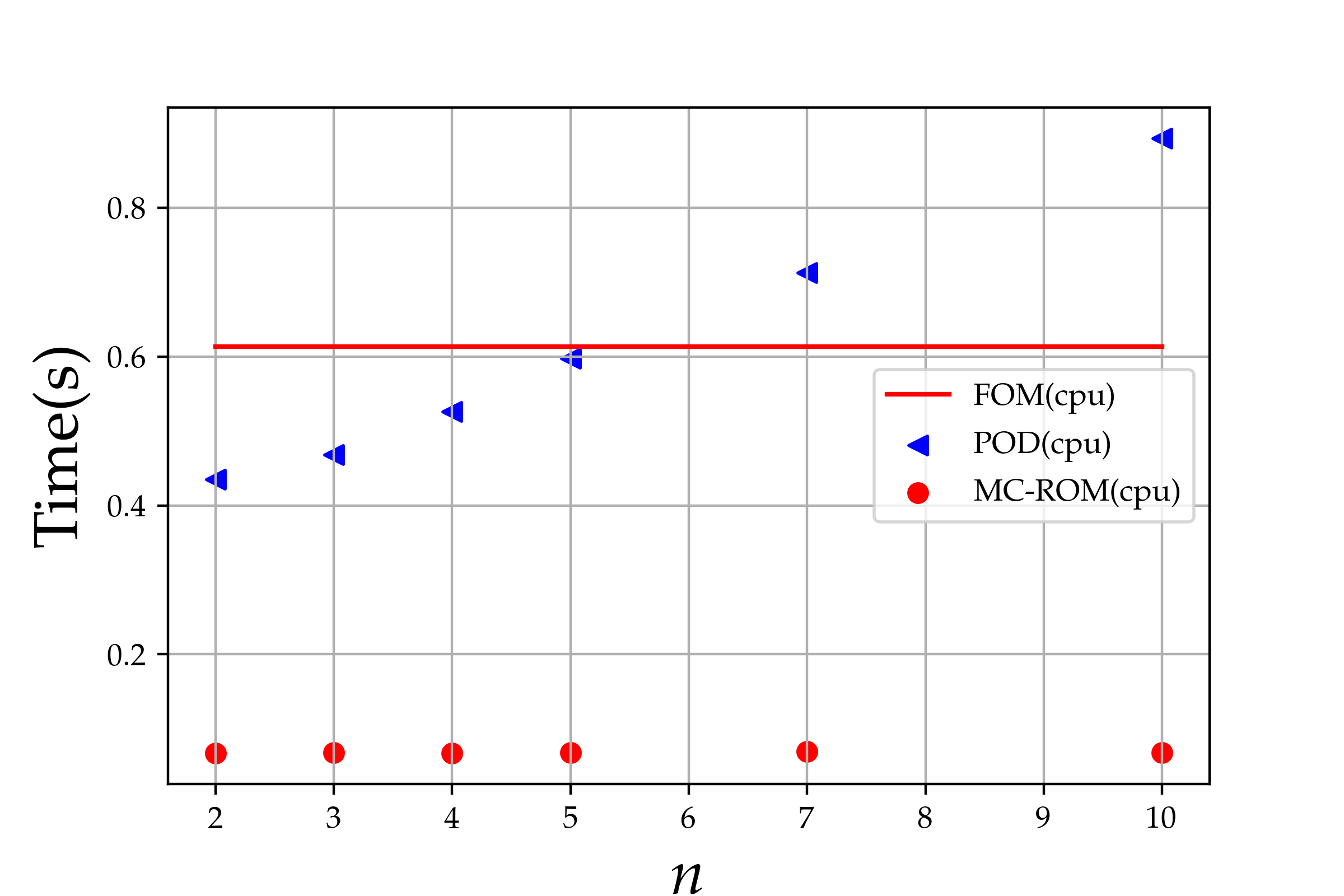}
    }
    \caption{Solving 1D viscous Burgers equations by the POD and the MC-ROM with viscosity term
    coefficient in $[0.5, 2]$. (a): Average relative error; (b): Average CPU time of
	online computation.}
    \label{fig:Burgers_diff}
\end{figure}

Figure \ref{subfig:a} shows the average relative error\,\eqref{total} of MC-ROM
and POD on the test set as a function of the dimension $n$ of the reduced space for
solving 1D viscous Burgers equations. The
POD method has better precision than the MC-ROM when $n \geq 4$. However, the MC-ROM
is still superior to the POD from two perspectives. One is that the online computational
performance of the MC-ROM is more efficient than that of the POD. Figure \ref{subfig:b} gives the
online inference time of the MC-ROM and the POD as a function of $n$. It can be seen that
the inference time of the MC-ROM is much lower than the POD and does not increase as the
dimension reduction space increases. However the inference time of the POD significantly
increases with a increase of $n$, and it will soon exceed FOM when $n > 6$. Another advantage is that the MC-ROM has a well generalization ability as the DL-ROM does for the
convection-dominated problems. As Figure \ref{fig:burgers_ad_test_error} shows, when the coefficient of the viscous
term decreases and the convection term becomes dominant, the POD method needs a more large $n (\geq 15)$  to obtain the same precision as the MC-ROM.

\subsection{2D parabolic equation}

\begin{figure}[h]
\centering
\includegraphics[width=0.4\textwidth]{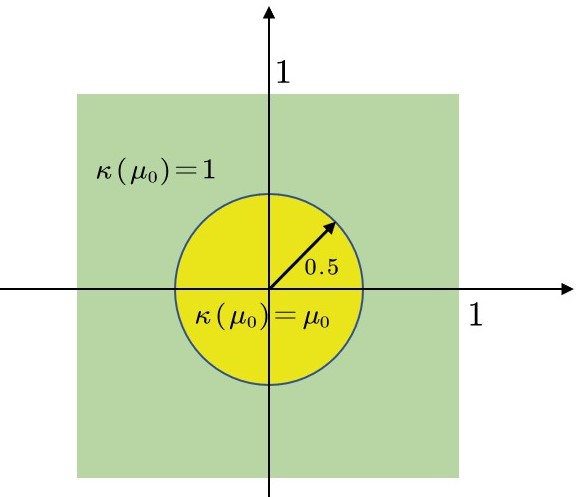}
\captionsetup{justification=centering}
\caption{Distribution of coefficient $\kappa$.}
\label{fig:kappa}
\end{figure}

Consider the following parabolic equation with discontinuous diffusion coefficient
\begin{equation}
     \left\{\begin{array}{ll}
     \partial_{t} u(t; \boldsymbol{\mu})-\operatorname{div}\left(\kappa\left(\mu_{0}\right) \nabla u(t; \boldsymbol{\mu} )\right)=0 & \text {in} \quad\Omega \times[0,T], \\
     u(t=0; \boldsymbol{\mu})=\mu_{1}(x-1)(x+1)(y-1)(y+1) & \text {in} \quad\Omega, \\
     u=0 \quad & \text {on} \quad\partial \Omega,
     \end{array}\right. 
     \label{parabolic}
\end{equation}
where $T=3$, $\Omega=(-1,1)^2$ is divided into $\Omega_1$ and $\Omega_2$
\begin{itemize}
     \item $\Omega_1$ is a disk centered at the origin of radius $r_0=0.5$, and
     \item $\Omega_2=\Omega/\ \overline{\Omega}_1$. 
 \end{itemize}
This equation contains two parameters $\mu_0$ and $\mu_1$. $\mu_0$ is related to the coefficient in $\Omega_1$, and $\mu_1$ appears in the initial value.
As shown in Figure \ref{fig:kappa}, the coefficient $\kappa$ is constant on $\Omega_1$ and $\Omega_2$, i.e.
\begin{equation}
\kappa|_{\Omega_1}=\mu_0 \quad \textrm{and} \quad \kappa|_{\Omega_2}=1.
\end{equation}

\begin{figure}[htp]
    \centering
    \includegraphics[width=0.4\textwidth]{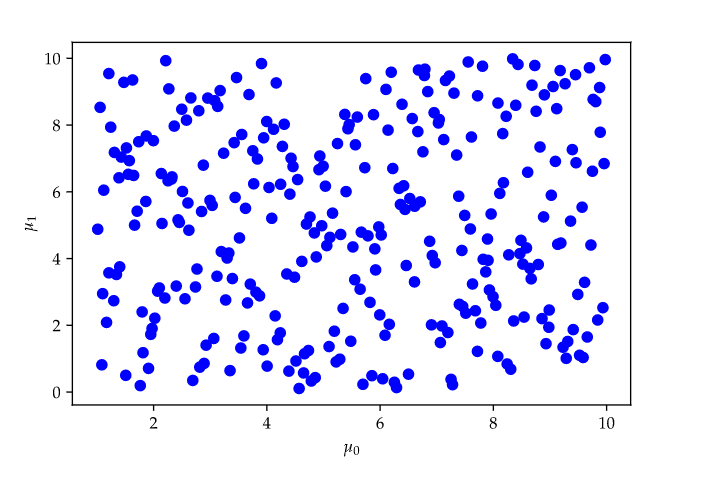}
    \hspace{0.5cm}
    \includegraphics[width=0.4\textwidth]{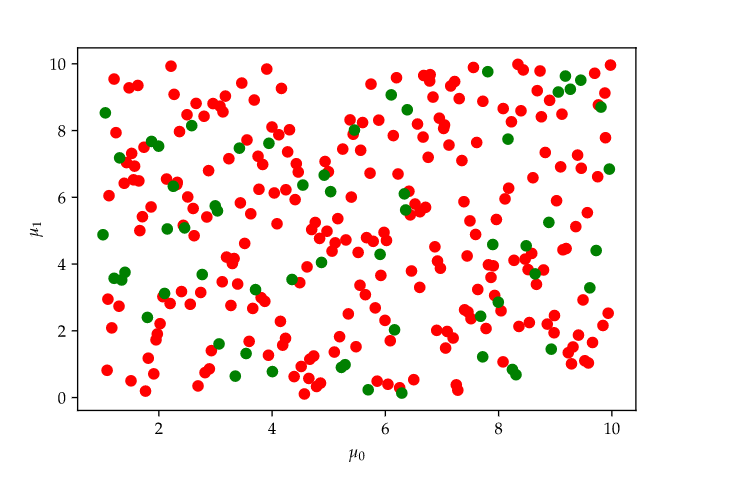}
    \caption{Left: Latin Hypercube Sampling with uniform distribution; Right: split left samples into training (red) and testing (green) samples.}
    \label{fig:sample}
\end{figure}

We consider parameter $\boldsymbol{\mu} = (\mu_0,\mu_1)\in \mathcal{P}=[1,10]\times[0.1,10].$ As Figure \ref{fig:sample} shows, we use
 Latin Hypercube Sampling \cite{stein1987large} with a uniform distribution to generate a total of $300$ parameters.

 According to the experience of the training process of the viscous Burgers equation, the MC-ROM is not overfitting. Therefore, we divide the dataset into a training set and test set with the ratio of 8:2, i.e., $240$ for training and $60$ for testing, no longer leave for the validation set. Then we solve the corresponding high-fidelity solution for training and testing. For a given parameter, we use FEM and backward Euler methods to discretize space and time variable, respectively. And we use the software \lstinline{FEniCS}\cite{alnaes2015fenics} to implement the above numerical methods and obtain the dataset. 
\begin{figure}[htp]
\centering
\subfigure[finite element mesh]{\label{fig:subfig:a}
      \includegraphics[width=0.27\textwidth]{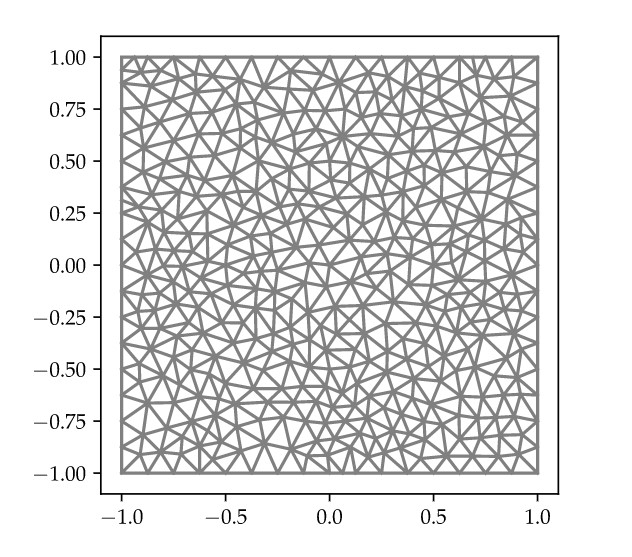}
}
\quad
\subfigure[$\boldsymbol{\mu} = (1.6472, 6.4912), t=0.0$]{
      \label{fig:subfig:b}
      \includegraphics[width=0.27\textwidth]{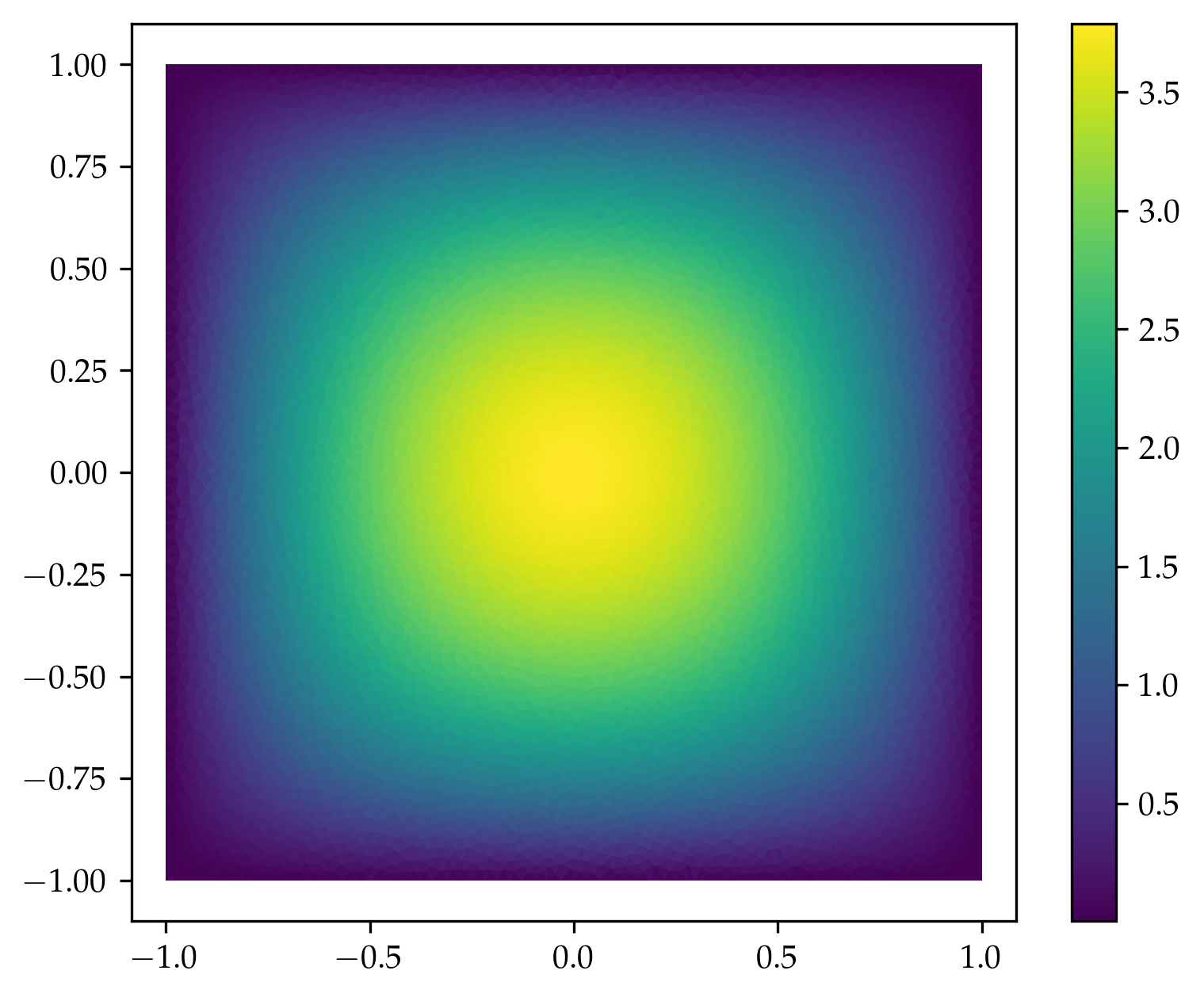}
}
\quad
\subfigure[$\boldsymbol{\mu} = (1.6472, 6.4912), t=3.0$]{
      \label{fig:subfig:c}
      \includegraphics[width=0.27\textwidth]{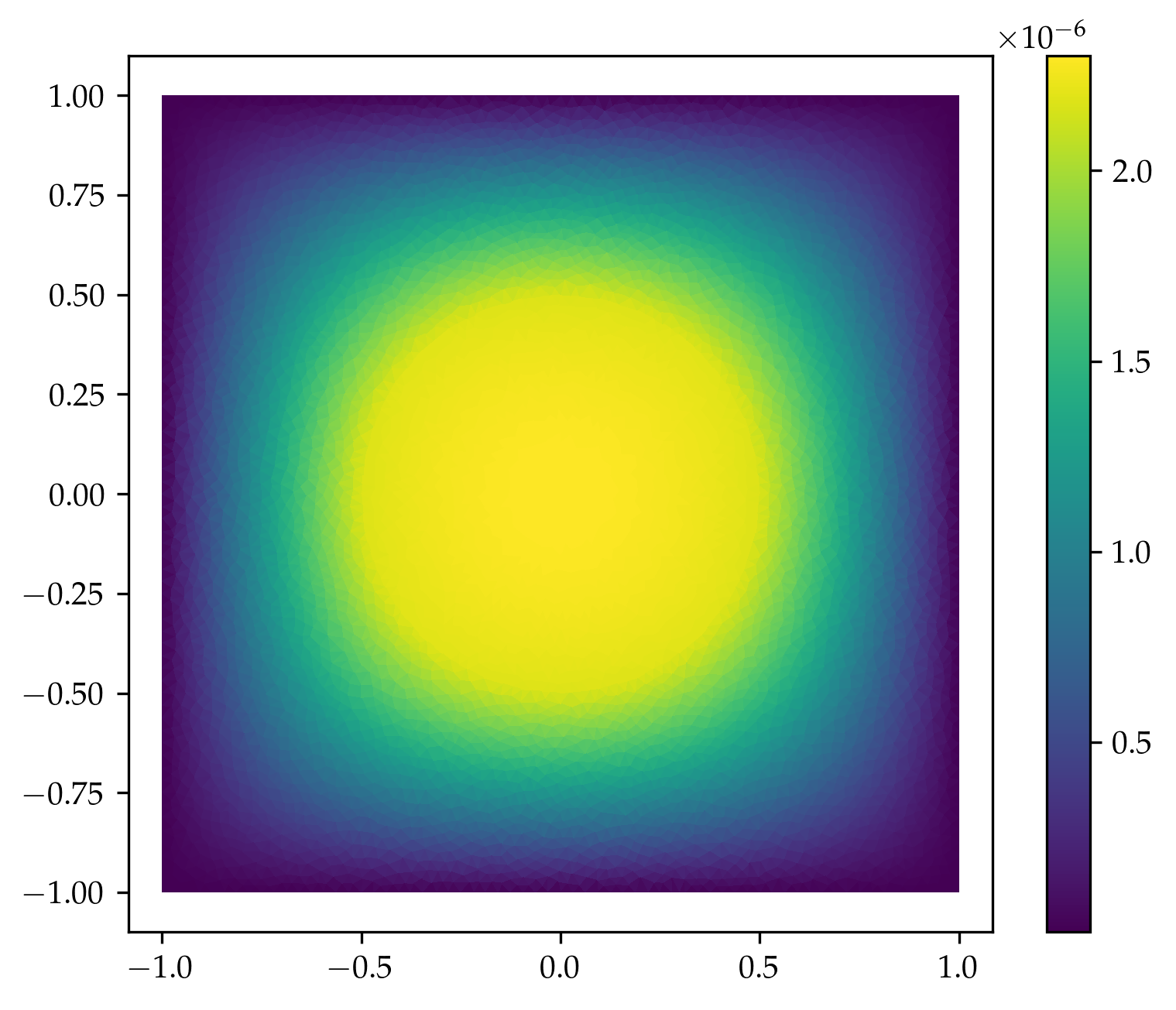}
}
\caption{(a): Schematic diagram of mesh; (b,c): Numerical solution of equation \eqref{parabolic} corresponding to $\boldsymbol{\mu} = (1.6472, 6.4912)$ at initial and final time.}
\label{fig:subfig}
\end{figure}
Figure \ref{fig:subfig:a} shows a schematic mesh of spatial discretization. In practice, the mesh we use contains $5716$ nodes and $11171$ cells and $\Delta t=0.05$. Figure \ref{fig:subfig:b} and Figure \ref{fig:subfig:c} show the FOM solution of parameter $(1.6472, 6.4912)$ at the initial and end moments. We can see that the order of magnitude of the solution is decreasing over time. Similar to the viscous Burgers equation, we classify training data using SVM. Table \ref{tab:parabolic_clf} presents the classification results.
\begin{table}[htp]
     \centering
     \caption{Classification results of the training set whose data are FOM solutions of the 2D parabolic equation.}
     \resizebox{\textwidth}{!}{
     \begin{tabular}{|l|c|c|c|c|c|c|}
     \hline
     Range of $\left\lVert\mathbf{u}_{h} \right\rVert_{\ell^{\infty}}$  &   $\left\lVert\mathbf{u}_{h}\right\rVert_{\ell^{\infty}}\geq 10^{-1}$  &   $[10^{-2}, 10^{-1})$     &  $[10^{-3}, 10^{-2})$       &  $[10^{-4}, 10^{-3})$       &   $[10^{-5}, 10^{-4})$      & $\left\lVert\mathbf{u}_{h}\right\rVert_{\ell^{\infty}}\leq 10^{-5}$ \\
     \hline
     Label  & 1      & 2      & 3      & 4      & 5      & 6 \\
     \hline
     Size     & 3751   & 2346 & 2368  & 2360  & 2360  & 1455\\
     \hline
     \end{tabular}}
    
     \label{tab:parabolic_clf}
 \end{table}%
 
The MC-ROM's subnets are all DL-ROM as shown in Figure \ref{fig:DL_ROM}, its concrete structure is as follows. $\Psi$ contains $10$ hidden layers, each layer contains $50$ neurons. Table \ref{tab:Network_structure_parabolic} shows the parameters of $\boldsymbol{h}_{\mathrm{enc}},\,\boldsymbol{h}_{\mathrm{dec}}$.
\begin{table}[htp] 
\centering
\caption{Network parameters of 2D Parabolic equation experiment. The meaning of parameters in Conv2d and ConvTranspose2d are: in channels(i), out channels(o), kernel size(k), stride(s), padding(p).}
\resizebox{\textwidth}{!}{
\begin{tabular}{|llr|ll|}
    \hline
    \multicolumn{2}{|c}{$\boldsymbol{h}_{\mathrm{enc}}$ (input shape:($N_b$, 5176))} &        & \multicolumn{2}{c|}{$\boldsymbol{h}_{\mathrm{dec}}$ (input shape:($N_b$, $n$))} \\
    \hline
    Layer type & Output shape &        & Layer type & Output shape \\
    \hline
    Fully connected &($N_b$, 4096)&        & Fully connected & ($N_b$,256)\\
    Reshape & ($N_b$,1,64,64) &        & Fully connected & ($N_b$,4096) \\
    Conv2d(i=1,o=8,k=5,s=1,p=2) & ($N_b$,8,64,64) &        & Reshape & ($N_b$,64,8,8) \\
    Conv2d(i=8,o=16,k=5,s=2,p=2) & ($N_b$,16,32,32) &        & ConvTranspose2d(i=64,o=64,k=5,s=3,p=2) & ($N_b$,64,22,22)  \\
    Conv2d(i=16,o=32,k=5,s=2,p=2) & ($N_b$,32,16,16) &        & ConvTranspose2d(i=64,o=32,k=5,s=3,p=2) & ($N_b$,32,64,64) \\
    Conv2d(i=32,o=64,k=5,s=2,p=2) & ($N_b$,64,8,8) &        & ConvTranspose2d(i=32,o=16,k=5,s=1,p=2) & ($N_b$,16,64,64) \\
    Reshape & ($N_b$,4096) &        &ConvTranspose2d(i=16,o=1,k=5,s=1,p=2) & ($N_b$,1,64,64)  \\
    Fully connected & ($N_b$,256) &        & Reshape & ($N_b$,4096) \\
    Fully connected & ($N_b$,$n$)  &        & Fully connected & ($N_b$, 5176) \\
    \hline
\end{tabular}}
\label{tab:Network_structure_parabolic}
\end{table}
We use the Adam algorithm, and the initial learning rate is set to $\eta = 0.001$. A multi-step decay learning rate adjustment strategy is used, i.e., the current learning rate is reduced by $0.1$ when epochs become $5000$ and $10000$. The maximum training epoch is $N_{\text{epochs}}=40000$, and the batch size is $N_b=5000$. We take $n=5$ and apply Algorithm \ref{algorithm:MC_ROM_offline} to train MC-ROM.

 Table \ref{fig:test_case} shows the solutions of PDE \eqref{parabolic} at time $t=0.2, 1.25, 2.20, 2.70$ when $\boldsymbol{\mu}=(9.9560,6.8453)$, obtained by FOM, MC-ROM, and DL-ROM, respectively. It also gives the relative error of MC-ROM and DL-ROM according to equation \eqref{single}. We find that MC-ROM maintains a good approximation at different times. At the same time, DL-ROM becomes worse over time, which means that MC-ROM has a better generalization ability than DL-ROM. Similar results can be found for other parameters $\boldsymbol{\mu}$ in the test set.
 \begin{table}[h]
    \centering
    \begin{tabular}{ | c | c | c | c|}
      \hline
       & FOM & MC-ROM & DL-ROM \\ \hline
     $ t = 0.2$& \begin{minipage}[b]{0.27\textwidth}
          \centering
          \raisebox{-.5\height}{\includegraphics[width=\linewidth]{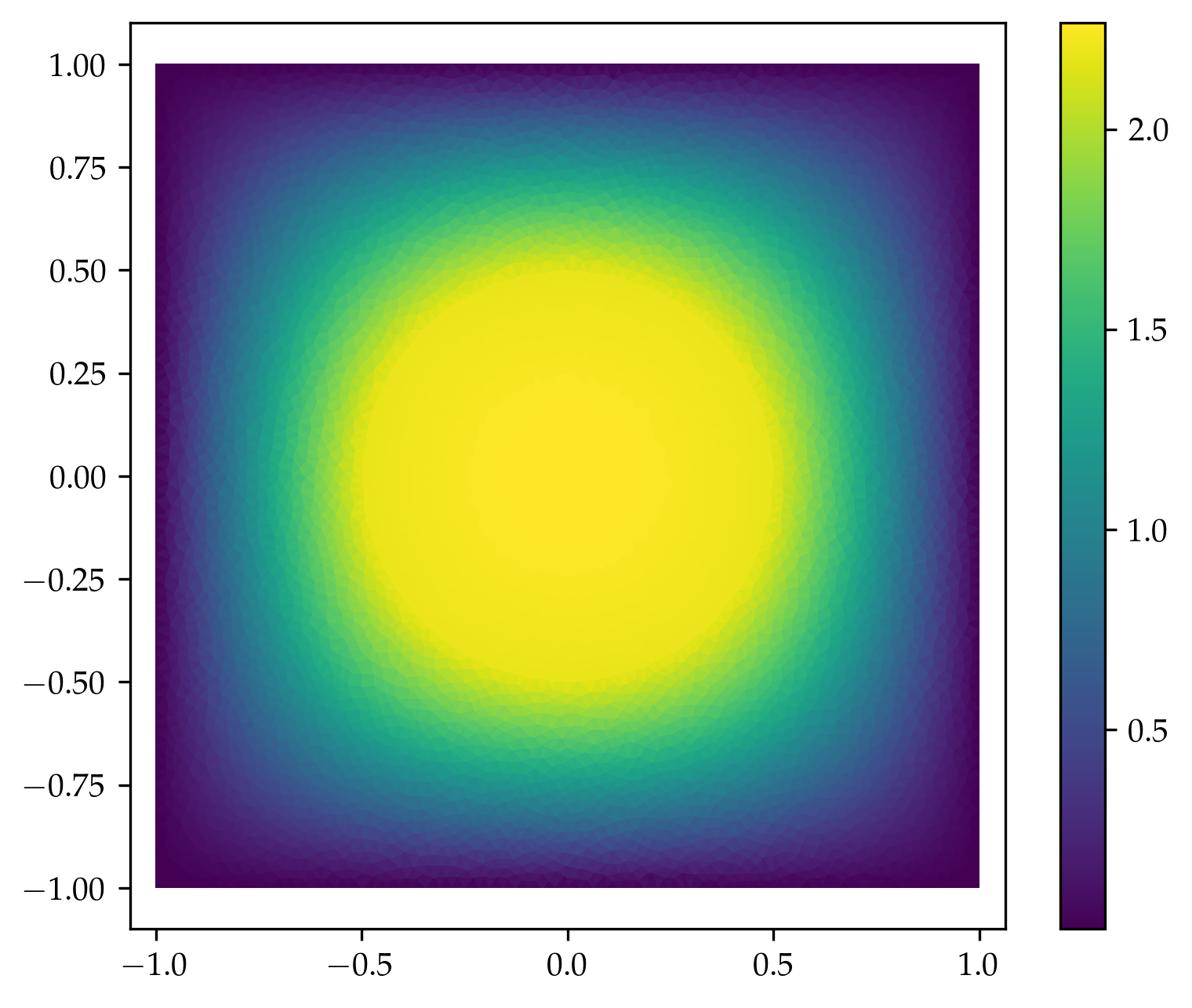}}
      \end{minipage}
      & \begin{minipage}[b]{0.27\textwidth}
        \centering
        \raisebox{-.5\height}{\includegraphics[width=\linewidth]{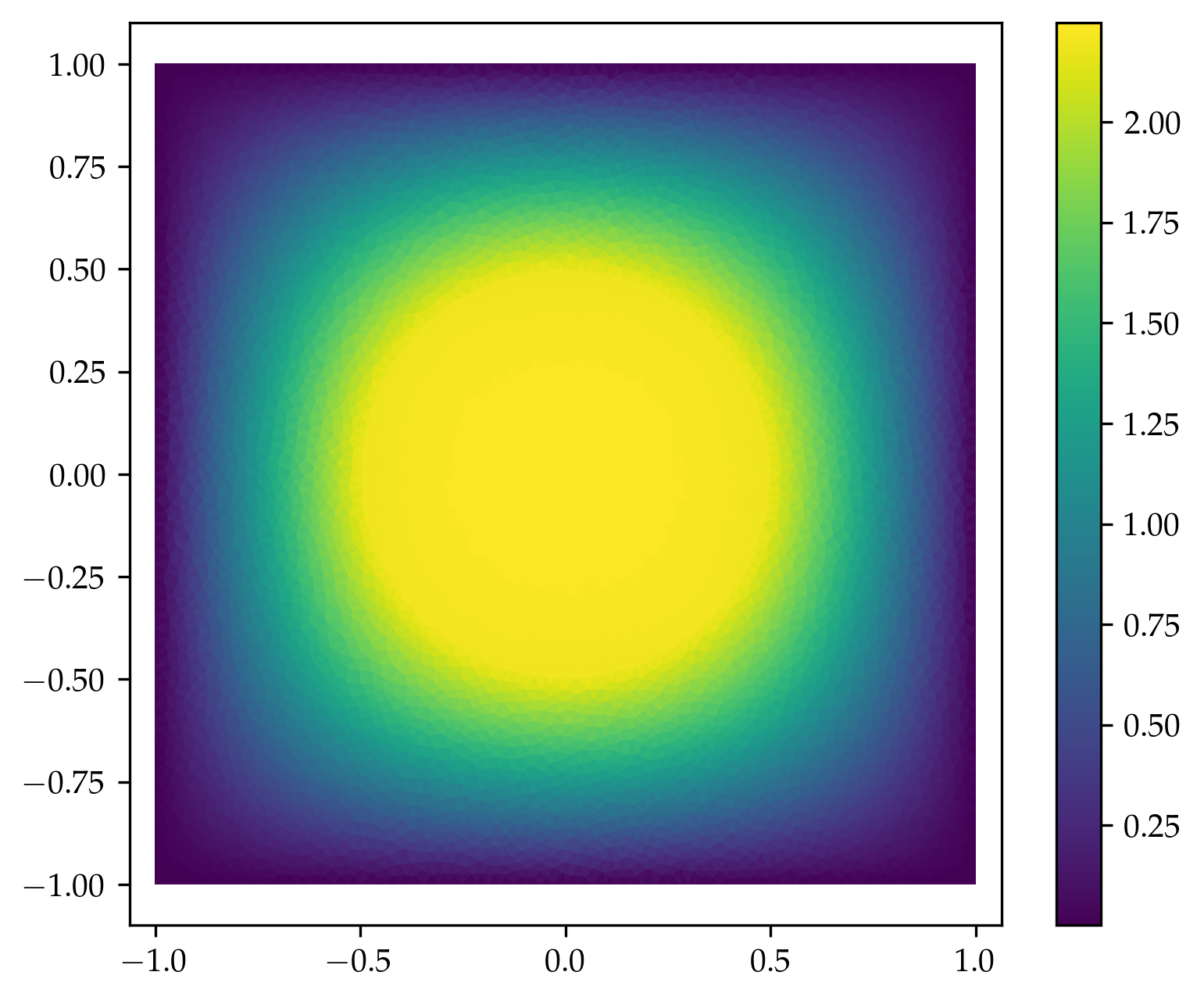}}
    \end{minipage}
      & \begin{minipage}[b]{0.27\textwidth}
        \centering
        \raisebox{-.5\height}{\includegraphics[width=\linewidth]{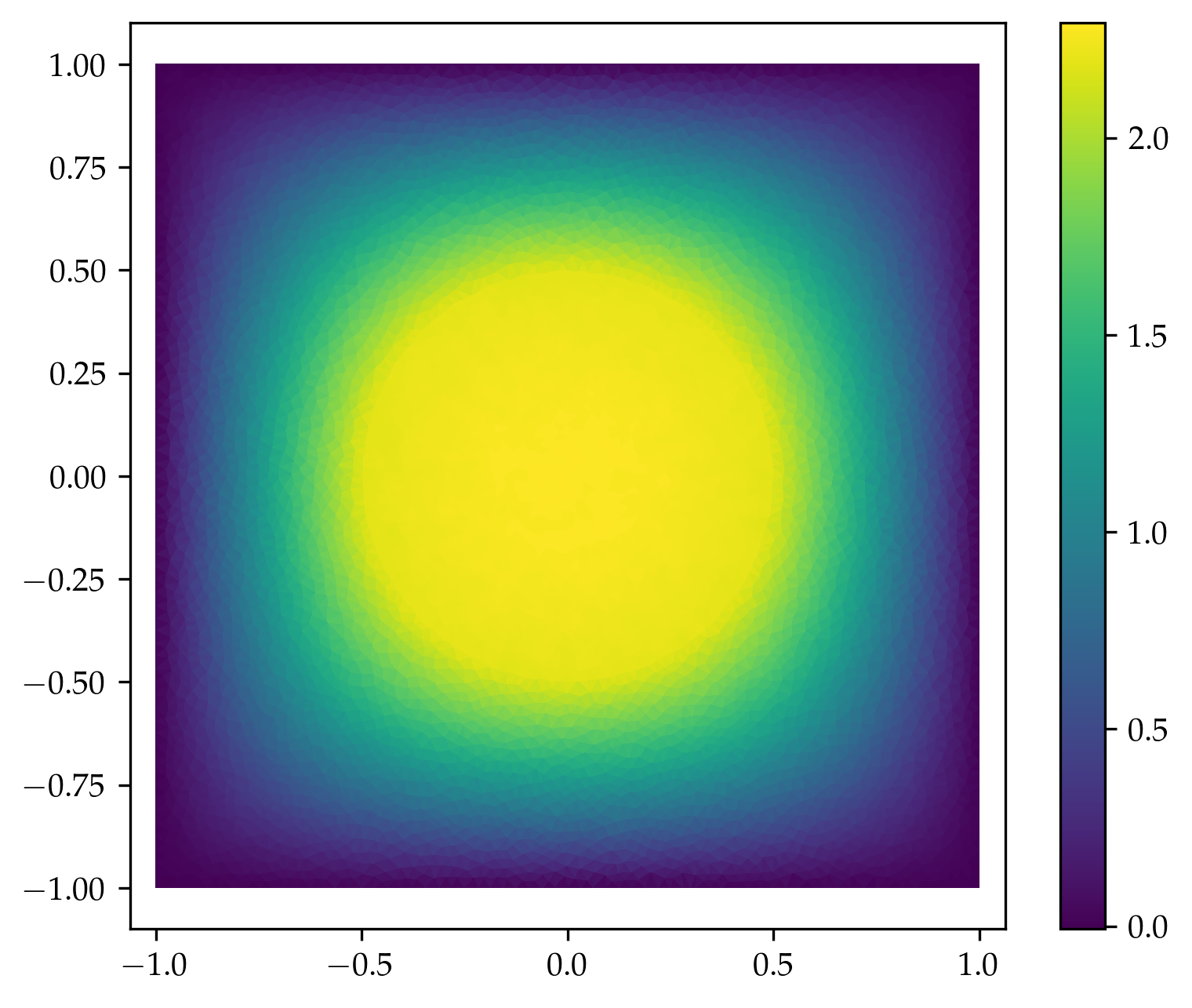}}
    \end{minipage}
      \\ \hline
    & & $\mathcal{E}_{single}=0.003609$ &  $\mathcal{E}_{single}=0.007461$   \\ \hline
      $ t = 1.25$& \begin{minipage}[b]{0.27\textwidth}
        \centering
        \raisebox{-.5\height}{\includegraphics[width=\linewidth]{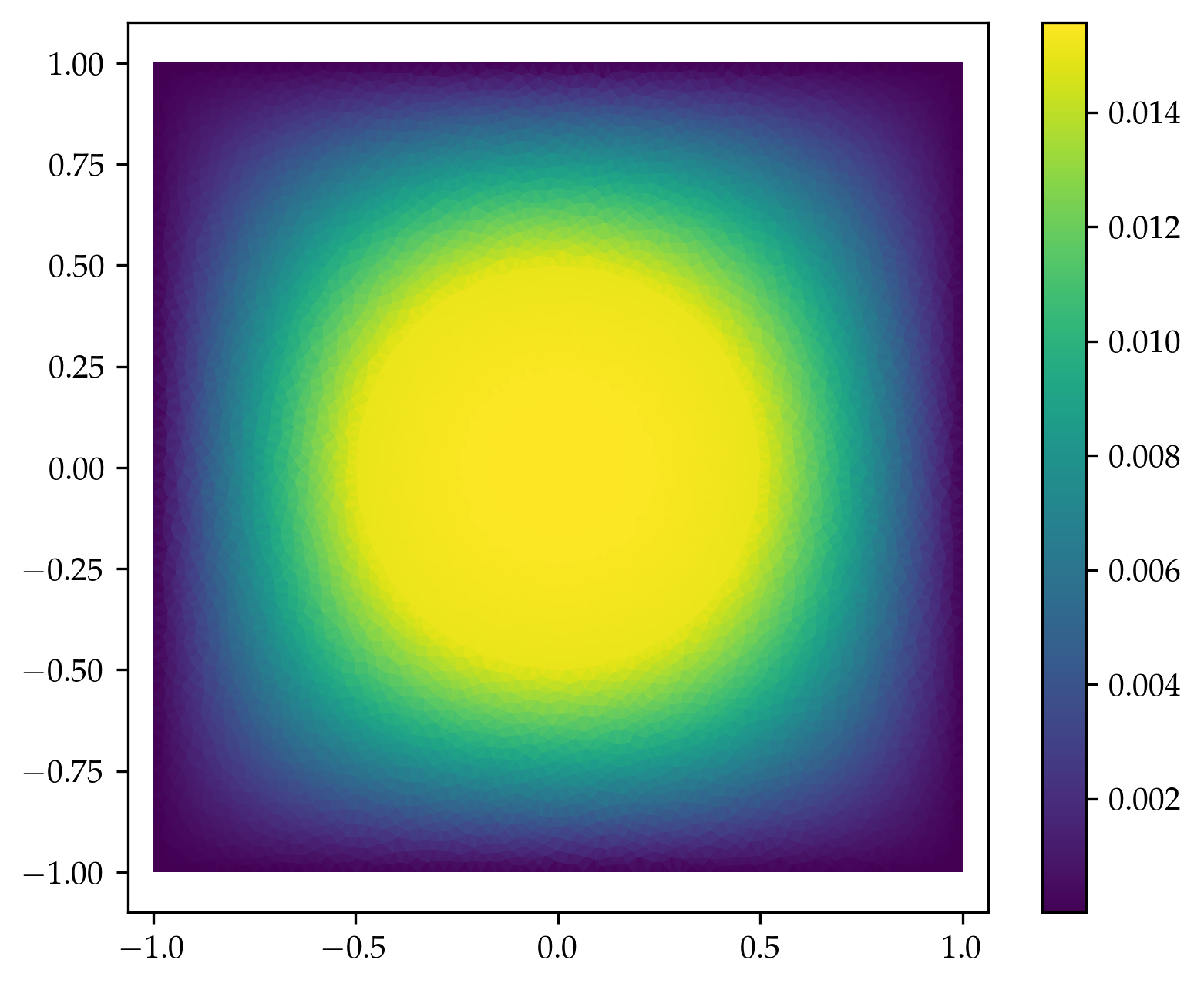}}
    \end{minipage}
    & \begin{minipage}[b]{0.27\textwidth}
      \centering
      \raisebox{-.5\height}{\includegraphics[width=\linewidth]{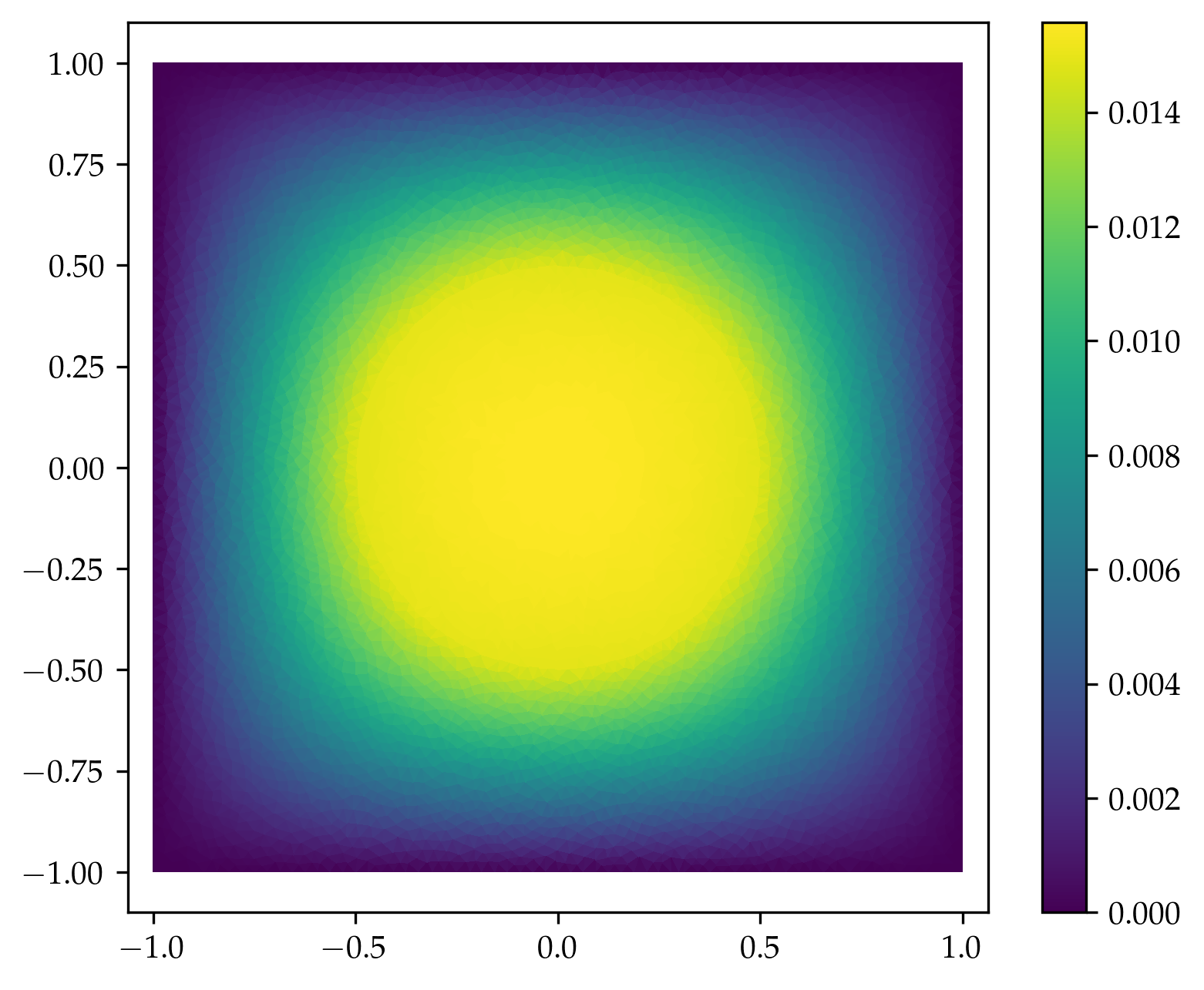}}
  \end{minipage}
    & \begin{minipage}[b]{0.27\textwidth}
      \centering
      \raisebox{-.5\height}{\includegraphics[width=\linewidth]{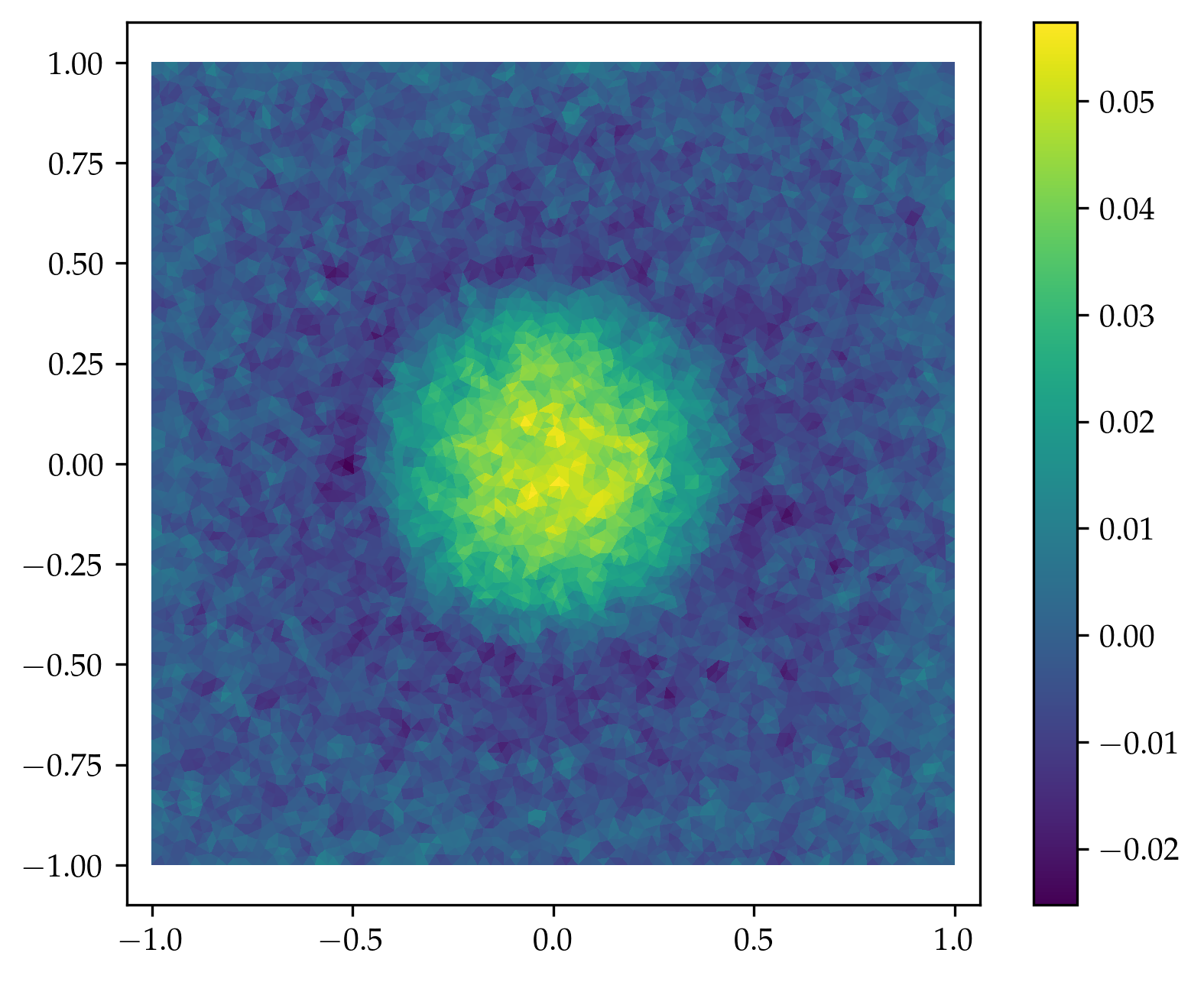}}
  \end{minipage}
    \\ \hline
    & & $\mathcal{E}_{single}=0.006973$ &  $\mathcal{E}_{single}=1.534366$   \\ \hline

    $ t = 2.20$& \begin{minipage}[b]{0.27\textwidth}
        \centering
        \raisebox{-.5\height}{\includegraphics[width=\linewidth]{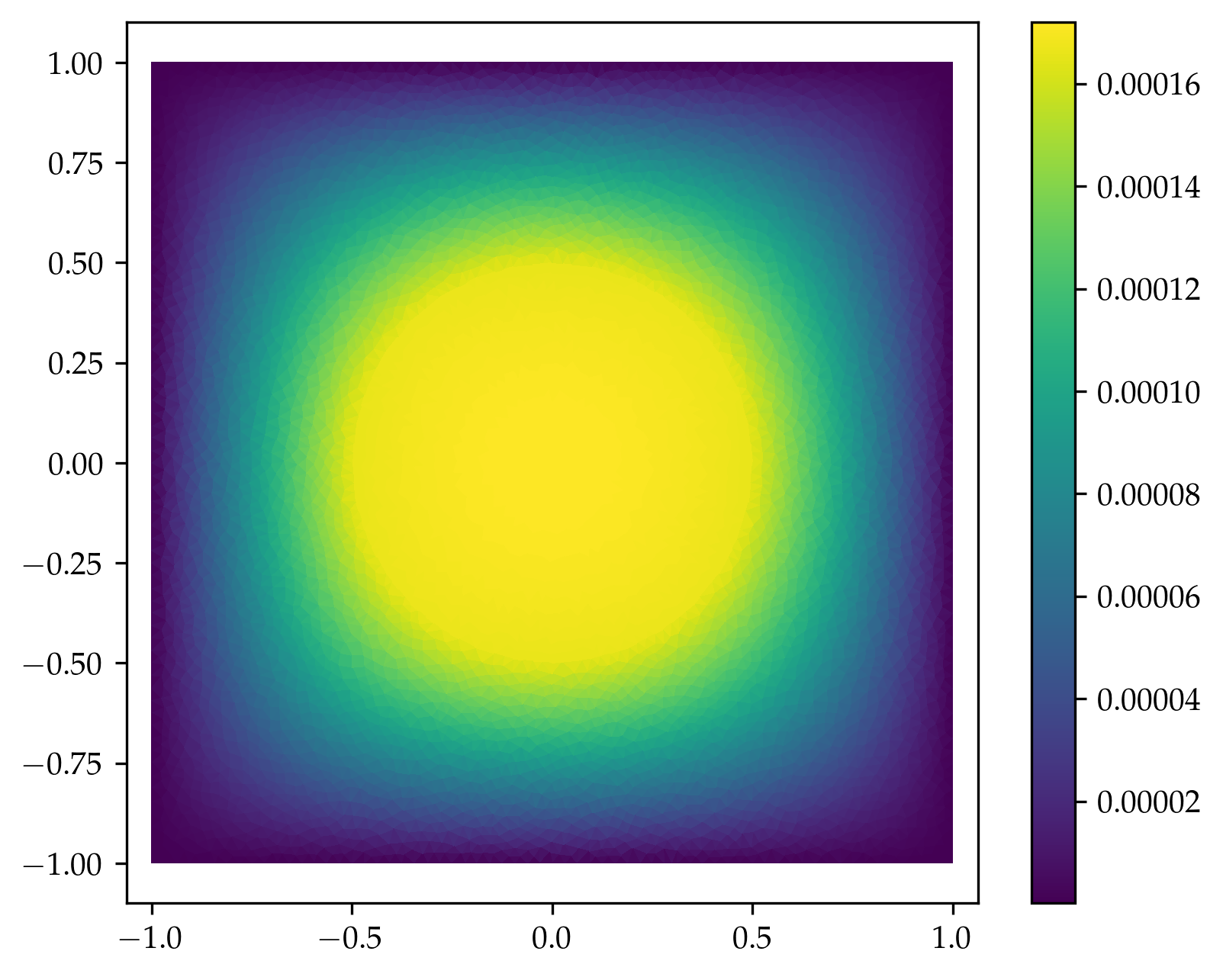}}
    \end{minipage}
    & \begin{minipage}[b]{0.27\textwidth}
      \centering
      \raisebox{-.5\height}{\includegraphics[width=\linewidth]{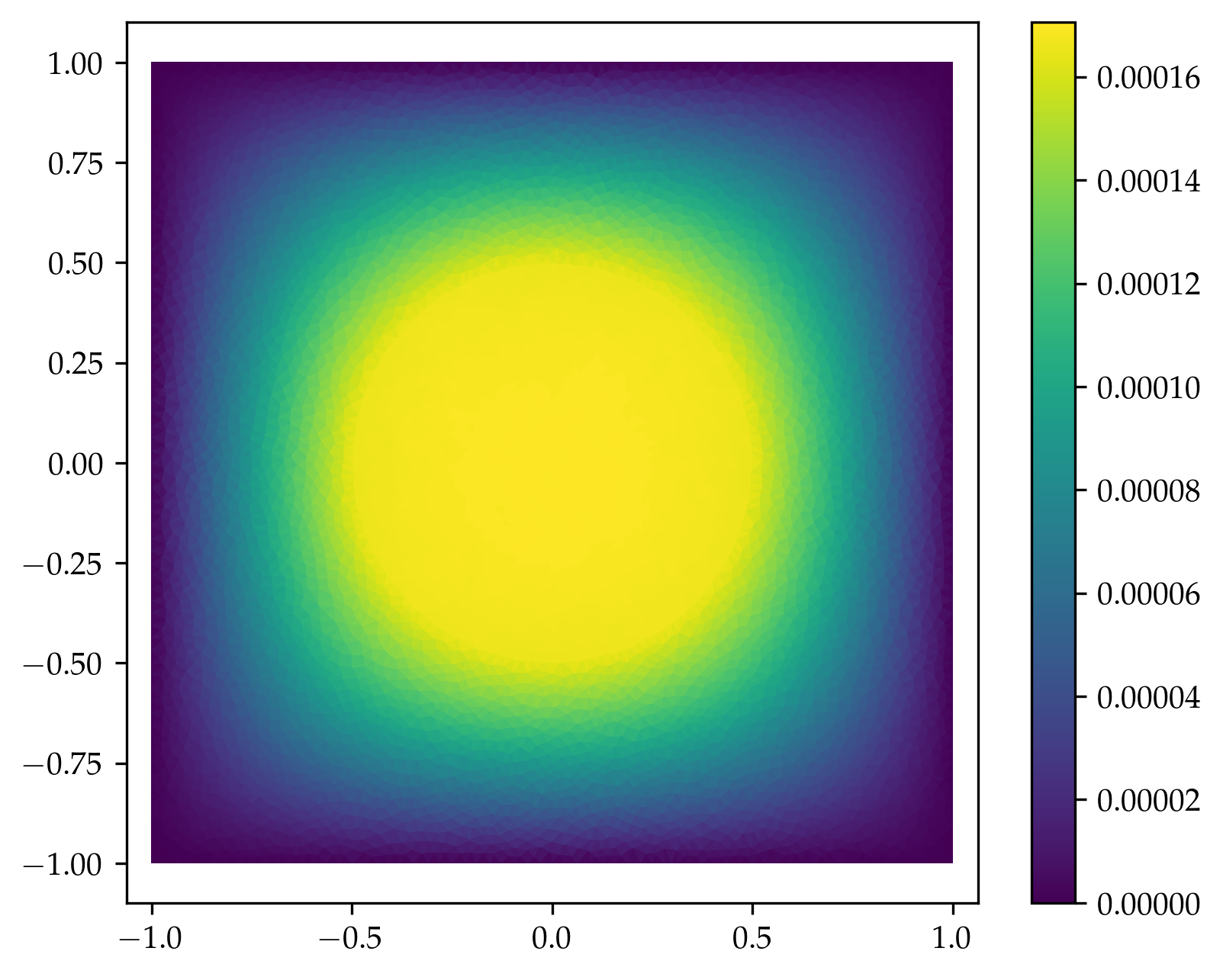}}
  \end{minipage}
    & \begin{minipage}[b]{0.27\textwidth}
      \centering
      \raisebox{-.5\height}{\includegraphics[width=\linewidth]{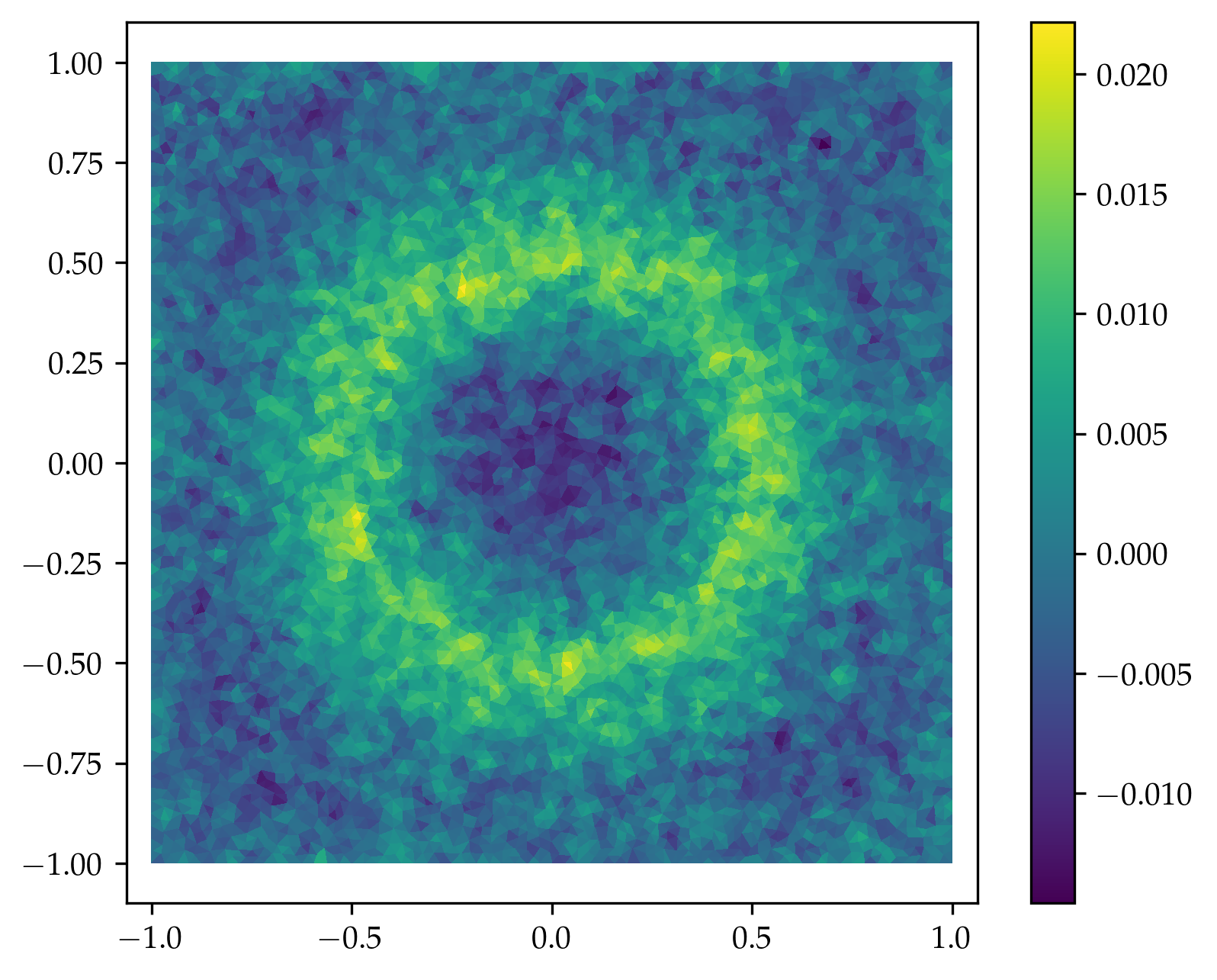}}
  \end{minipage}
    \\ \hline
    & & $\mathcal{E}_{single}=0.003970$ &  $\mathcal{E}_{single}=63.122426$   \\ \hline

    $ t = 2.70$& \begin{minipage}[b]{0.27\textwidth}
        \centering
        \raisebox{-.5\height}{\includegraphics[width=\linewidth]{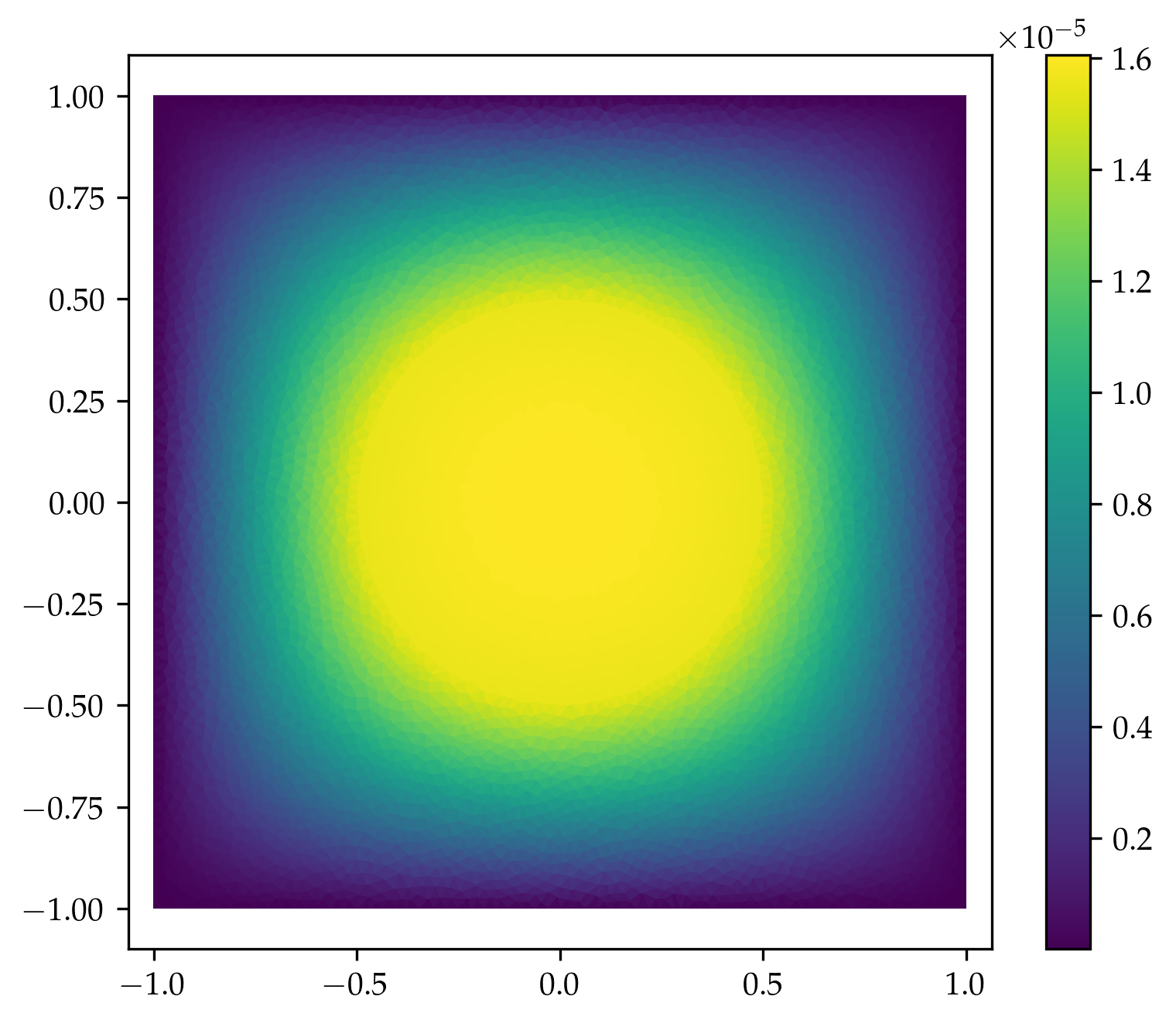}}
    \end{minipage}
    & \begin{minipage}[b]{0.27\textwidth}
      \centering
      \raisebox{-.5\height}{\includegraphics[width=\linewidth]{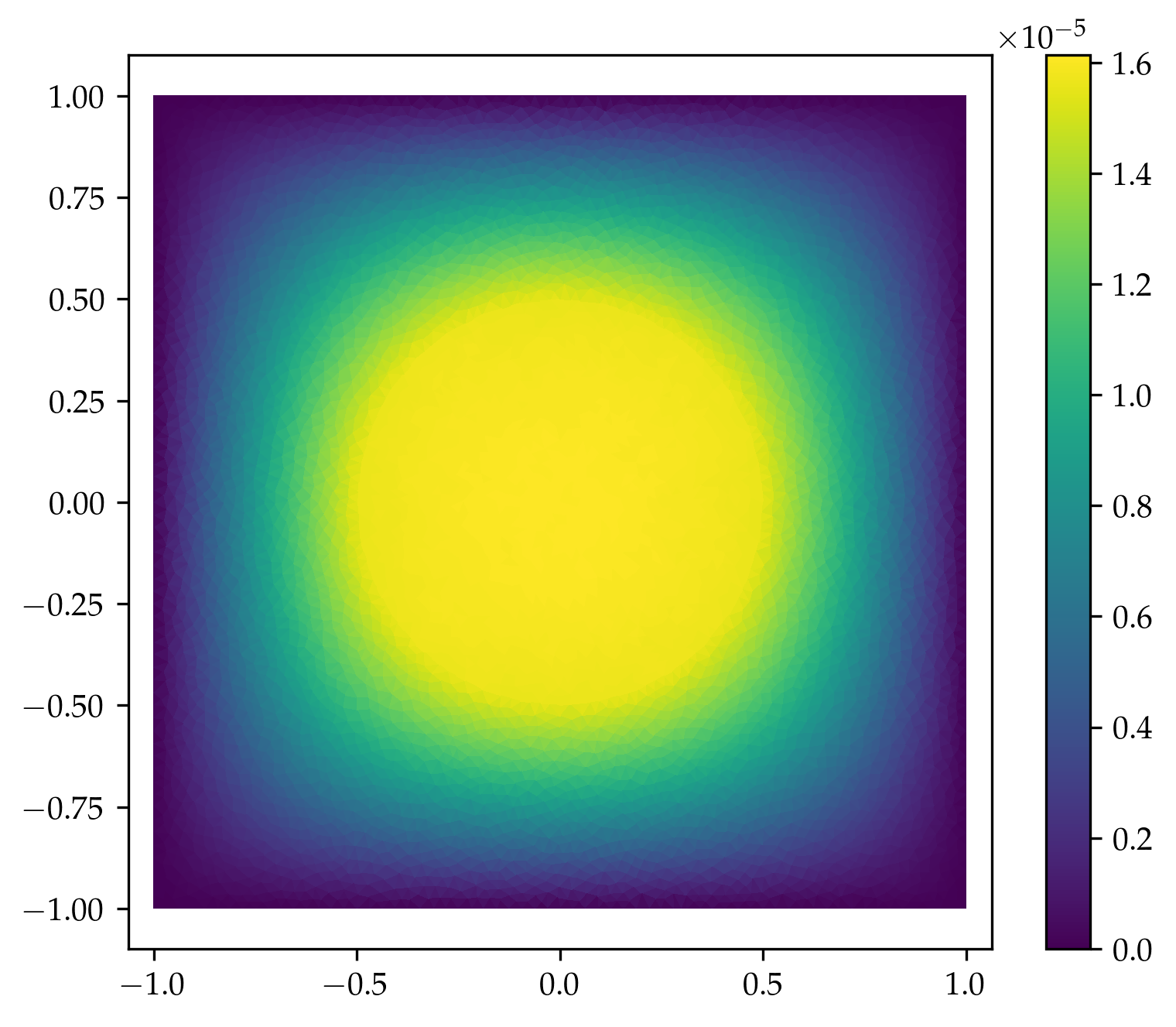}}
  \end{minipage}
    & \begin{minipage}[b]{0.27\textwidth}
      \centering
      \raisebox{-.5\height}{\includegraphics[width=\linewidth]{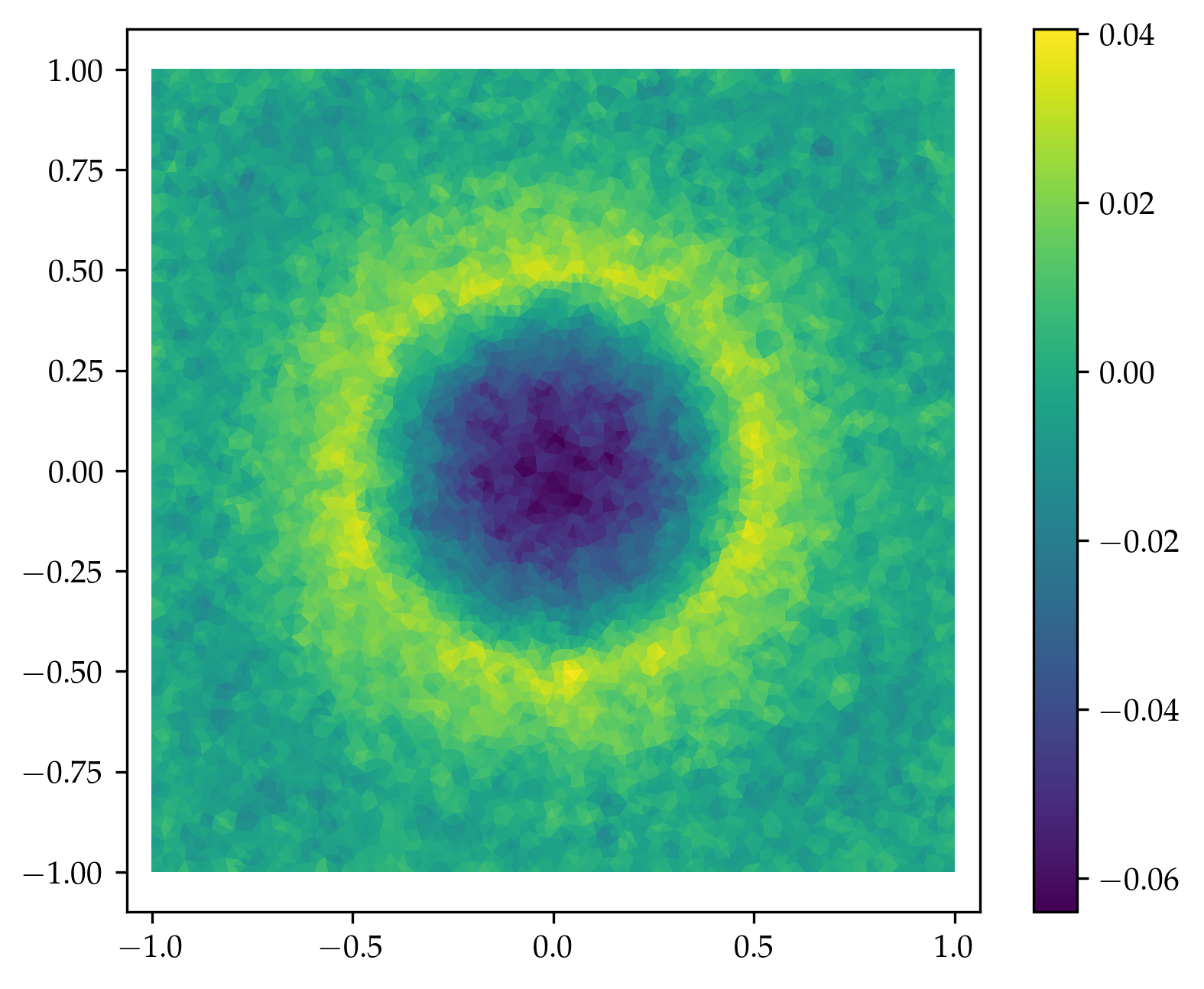}}
  \end{minipage}
    \\ \hline
    & & $\mathcal{E}_{single}=0.008112$ &  $\mathcal{E}_{single}=1672.197287$   \\ \hline

    \end{tabular}
    \caption{Numerical results of solving PDE\,\eqref{parabolic} when $\boldsymbol{\mu}=(9.9560,6.8453)$ by the FOM, the MC-ROM, and the DL-ROM at different times. The relative error  $\mathcal{E}_{single}$ is defined in \eqref{single}.}
    \label{fig:test_case}
  \end{table}

\begin{figure}[htp]
    \centering
    \subfigure[test error]{
        \label{fig:thermal_block_pod:a}
        \includegraphics[width=0.45\textwidth]{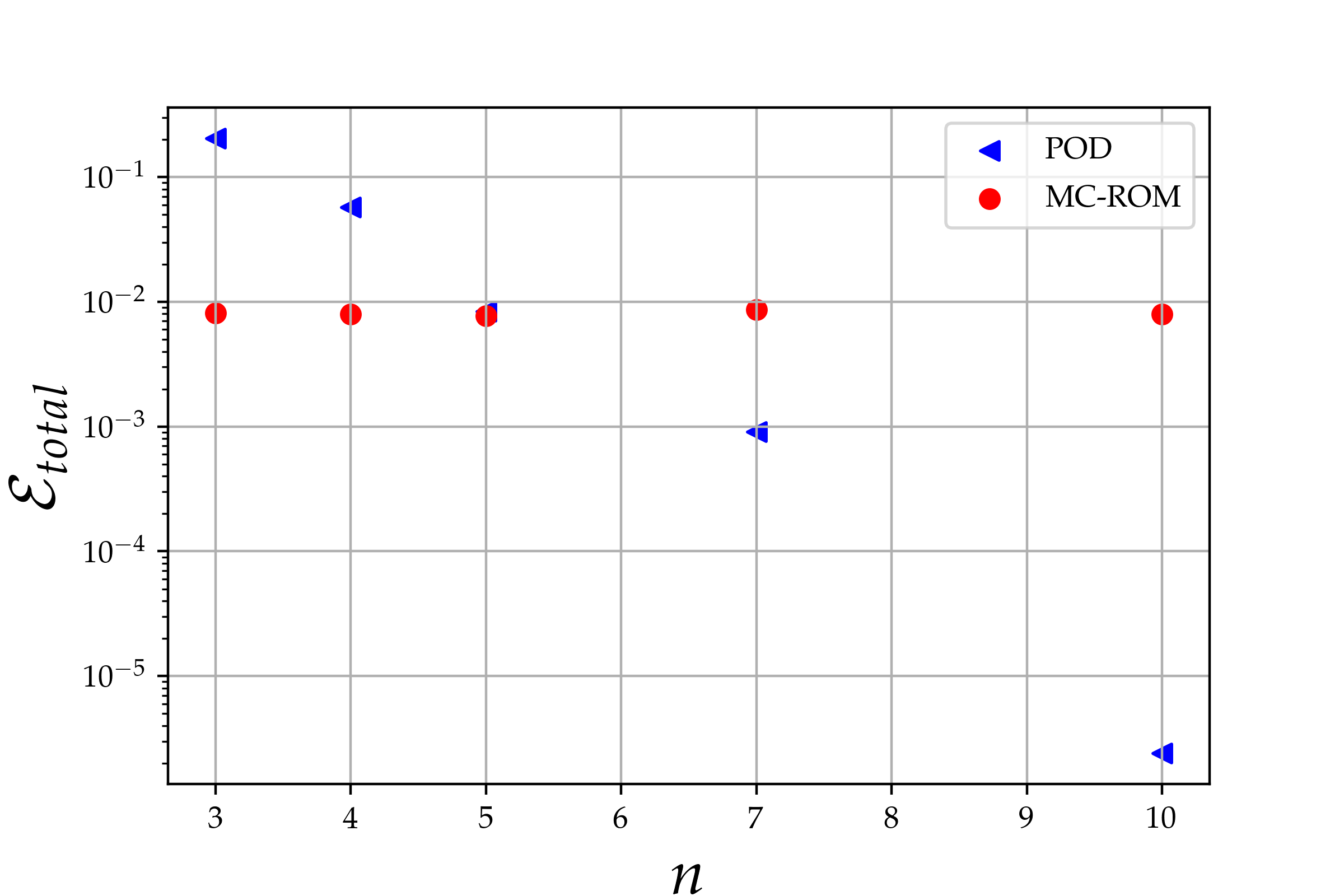}
    }
    \quad
    \subfigure[inference time]{
            \label{fig:thermal_block_pod:b}
            \includegraphics[width=0.45\textwidth]{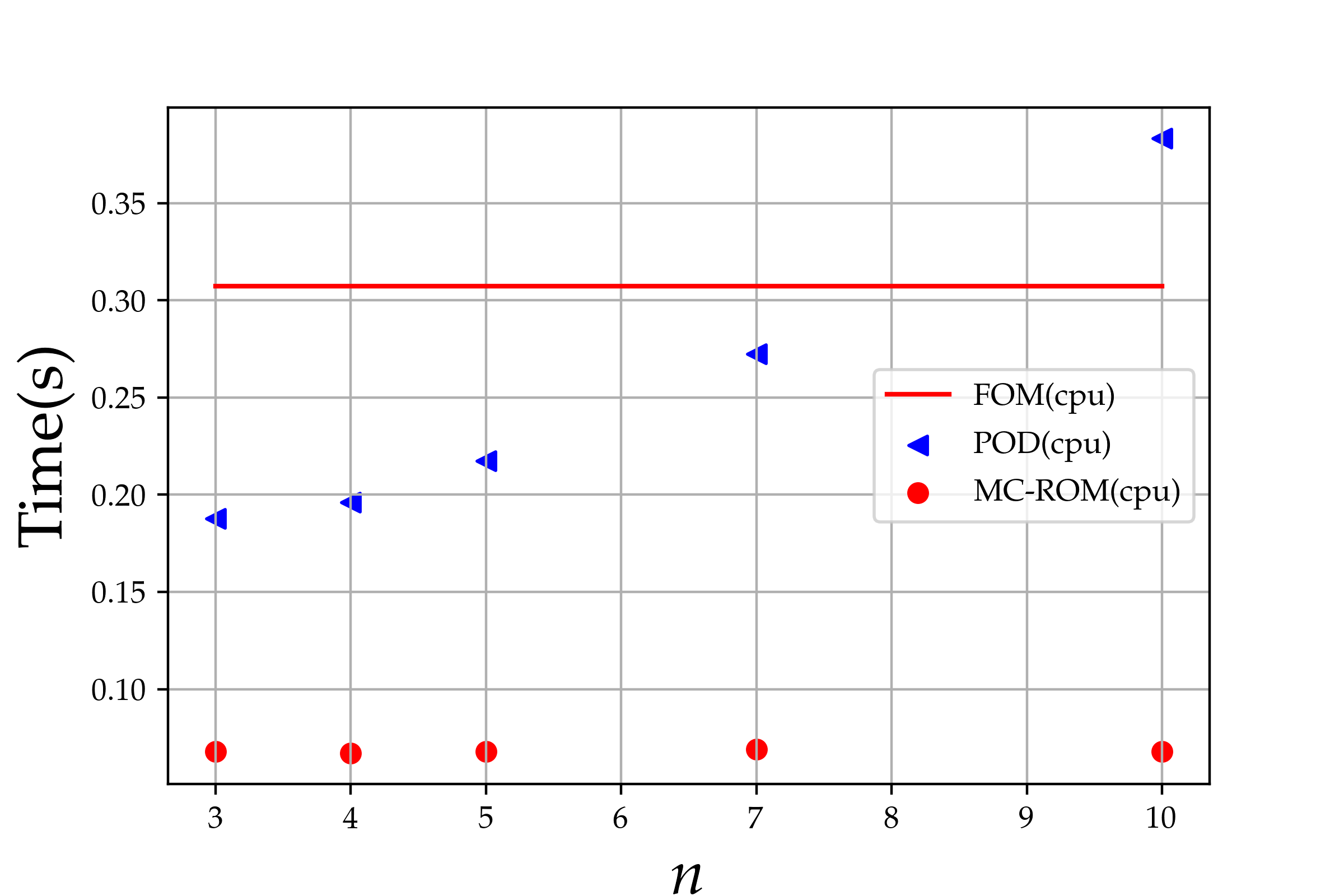}
    }
    \caption{Solving 2D parabolic equations by the POD and the MC-ROM (a): Average relative
	error; (b): Average CPU time of online computation.}
    \label{fig:thermal_block_pod}
\end{figure}

We also compared the POD and the MC-ROM for solving 2D parabolic equations, as shown
in Figure \ref{fig:thermal_block_pod}. Figure \ref{fig:thermal_block_pod:a} shows the
average relative error of the MC-ROM and the POD on the test set as a function of the
dimension $n$ of the reduced space. We can find the similar phenomena as the above 1D
viscous Burgers equations. The MC-ROM has better approximation than the POD when
$n\leq 5$, and the less online inference time for any $n$.

\section{Conclusion and discussions}\label{sec07}

This paper proposes the MC-ROM for solving time-dependent PPDEs. MC-ROM classifies
the data according to the magnitude of the FOM solutions, constructs corresponding
subnets for different data classes, and builds a classifier to integrate all subnets.
In the offline stage, MC-ROM can use transfer learning techniques to accelerate the
training of subnets with the same architecture. Numerical experiments show that
MC-ROM has good generalization ability for both diffusion-dominant and
convection-dominant problems. Compared with DL-ROM, MC-ROM maintains its advantages,
including its good approximation and  generalization ability for convection-dominant
problems, efficient online calculation. More importantly, MC-ROM improves the
generalization ability of DL-ROM for diffusion-dominant problems. And MC-ROM is much
more efficient than POD when doing online computing. Even more, POD is not suitable
for convection-dominant problems due to its poor approximation, while MC-ROM can
maintain sufficient accuracy.  In this work, we use the DL-ROM as subset due to  its
good approximation ability and efficient online calculation.  In fact, the MC-ROM is
a flexibly extended model.  Any appropriate network that approximates the
parameter-to-solution mapping can be used as a subnet in the MC-ROM.

There are still many works to be explored. For instance, the characteristic of PPDEs
that this paper focuses on is that the magnitude of their solutions changes
drastically along the time direction. Therefore, we regard time variable as a parameter that
can solve this kind of PPDEs. However, there are many PPDEs whose solutions change
drastically in space or in time-space, such as interface problems and multi-scale
problems. In that case, we need to regard spatial variables as parameters. 
To train the network to achieve satisfactory approximation accuracy, the
amount of required data increases exponentially as the dimension of parameter
space increases which leads to the course of dimension. 
Developing an efficient sampling approach to high-dimensional parameter space needs
to be solved in our further work.
Furthermore, we will apply MC-ROM to more scientific and engineering problems.

\bibliographystyle{unsrt}
\bibliography{references}

\end{document}